\documentclass[titlepage,draft,12pt]{article} 
\usepackage{pstricks,pst-plot,pst-node} 
\usepackage{amssymb,amsthm} 
\usepackage[a4paper]{geometry}

\date{}

\newcommand{\ep}{\varepsilon}
\renewcommand{\qed}{{\penalty 10000\mbox{$\quad\Box$}}}
\newcommand{\re}{\mathbb{R}}

\newcommand{\n}{\mathbb{N}}
\newcommand{\cep}{c_{\ep}}

\newcommand{\Eep}{E_{\ep}}

\newcommand{\lk}{\lambda_{k}}
\newcommand{\lku}{\lambda_{k-1}}
\newcommand{\oa}{a_{1}}
\newcommand{\ob}{b_{1}}
\newcommand{\ul}{u_{\lambda}}
\newcommand{\G}{\mathcal{G}}

\newtheorem{thm}{Theorem}[section]
\newtheorem{thmbibl}{Theorem}

\newtheorem{rmk}[thm]{Remark}

\newtheorem{defn}[thm]{Definition}

\newtheorem{lemma}[thm]{Lemma}

\title{Linear hyperbolic equations with time-dependent propagation
speed and strong damping}

\author{Marina Ghisi\vspace{1ex}\\ 
{\normalsize Universit\`a degli Studi di Pisa} \\
{\normalsize Dipartimento di Matematica}\\ 
{\normalsize PISA (Italy)}\\
{\normalsize e-mail: \texttt{ghisi@dm.unipi.it}}
\and
Massimo Gobbino\vspace{1ex}\\ 
{\normalsize Universit\`a degli Studi di Pisa} \\
{\normalsize Dipartimento di Matematica}\\ 
{\normalsize PISA (Italy)}\\  
{\normalsize e-mail: \texttt{m.gobbino@dma.unipi.it}}
}

\begin{document}
\maketitle
\begin{abstract}
	We consider a second order linear equation with a time-dependent
	coefficient $c(t)$ in front of the ``elastic'' operator.  For
	these equations it is well-known that a higher space-regularity of
	initial data compensates a lower time-regularity of $c(t)$.
	
	In this paper we investigate the influence of a strong
	dissipation, namely a friction term which depends on a power of
	the elastic operator.
	
	What we discover is a threshold effect.  When the exponent of the
	elastic operator in the friction term is greater than 1/2, the
	damping prevails and the equation behaves as if the
	coefficient $c(t)$ were constant.  When the exponent is less than
	1/2, the time-regularity of $c(t)$ comes into play.  If $c(t)$ is
	regular enough, once again the damping prevails.  On the contrary,
	when $c(t)$ is not regular enough the damping might be
	ineffective, and there are examples in which the dissipative
	equation behaves as the non-dissipative one.  As expected, the
	stronger is the damping, the lower is the time-regularity
	threshold.
	
	We also provide counterexamples showing the optimality of our
	results.
	
\vspace{6ex}

\noindent{\bf Mathematics Subject Classification 2010 (MSC2010):}
35L20, 35L80, 35L90.


\vspace{6ex}


\end{abstract}

 
\section{Introduction}

Let $H$ be a separable real Hilbert space.  For every $x$ and $y$ in
$H$, $|x|$ denotes the norm of $x$, and $\langle x,y\rangle$ denotes
the scalar product of $x$ and $y$.  Let $A$ be a self-adjoint linear
operator on $H$ with dense domain $D(A)$.  We assume that $A$ is
nonnegative, namely $\langle Ax,x\rangle\geq 0$ for every $x\in D(A)$,
so that for every $\alpha\geq 0$ the power $A^{\alpha}x$ is defined
provided that $x$ lies in a suitable domain $D(A^{\alpha})$.

We consider the second order linear evolution equation
\begin{equation}
	u''(t)+2\delta A^{\sigma}u'(t)+c(t)Au(t)=0,
	\label{pbm:eqn}
\end{equation}
with initial data
\begin{equation}
	u(0)=u_{0},
	\hspace{3em}
	u'(0)=u_{1}.
	\label{pbm:data}
\end{equation}

As far as we know, this equation has been considered in the literature
either in the case where $\delta=0$, or in the case where
$\delta>0$ but the coefficient $c(t)$ is constant.  Let us give a
brief outline of the previous literature which is closely related to
our results.

\paragraph{\textmd{\emph{The non-dissipative case}}}

When $\delta=0$, equation (\ref{pbm:eqn}) reduces to
\begin{equation}
	u''(t)+c(t)Au(t)=0.
	\label{pbm:non-diss}
\end{equation}

This is the abstract setting of a wave equation in which $c(t)$ 
represents the square of the propagation speed.

If the coefficient $c(t)$ is Lipschitz continuous and satisfies the
strict hyperbolicity condition
\begin{equation}
	0<\mu_{1}\leq c(t)\leq\mu_{2},
	\label{hp:sh}
\end{equation}
then it is well-know that problem
(\ref{pbm:non-diss})--(\ref{pbm:data}) is well-posed in the classic
energy space $D(A^{1/2})\times H$ (see for example the classic 
reference~\cite{LM}). 

If the coefficient is not Lipschitz continuous, things are more
complex, even if (\ref{hp:sh}) still holds true.  This problem was
addressed by F.~Colombini, E.~De~Giorgi and S.~Spagnolo in the seminal
paper~\cite{dgcs}.  Their results can be summed up as follows (we
refer to section~\ref{sec:notation} below for the precise
functional setting and rigorous statements).
\begin{enumerate}
	\renewcommand{\labelenumi}{(\arabic{enumi})} 
	
	\item Problem (\ref{pbm:non-diss})--(\ref{pbm:data}) has always a
	unique solution, up to admitting that this solution takes its
	values in a very large Hilbert space (ultradistributions).  This
	is true for initial data in the energy space $D(A^{1/2})\times H$,
	but also for less regular data, such as distributions or
	ultradistributions.

	\item If initial data are regular enough, then the solution is
	regular as well.  How much regularity is required depends on the
	time-regularity of $c(t)$.  Classic examples are the
	following.  If $c(t)$ is just measurable, problem
	(\ref{pbm:non-diss})--(\ref{pbm:data}) is well-posed in the class
	of analytic functions.  If $c(t)$ is $\alpha$-H\"{o}lder
	continuous for some $\alpha\in(0,1)$, problem
	(\ref{pbm:non-diss})--(\ref{pbm:data}) is well-posed in the Gevrey
	space of order $(1-\alpha)^{-1}$.

	\item If initial data are not regular enough, then the solution
	may exhibit a severe derivative loss for all positive times.  For
	example, for every $\alpha\in(0,1)$ there exist a coefficient
	$c(t)$ which is $\alpha$-H\"{o}lder continuous, and initial data
	$(u_{0},u_{1})$ which are in the Gevrey class of order $\beta$ for
	every $\beta>(1-\alpha)^{-1}$, such that the corresponding
	solution to (\ref{pbm:non-diss})--(\ref{pbm:data}) (which exists
	in the weak sense of point~(1)) is not even a distribution for
	every $t>0$.
\end{enumerate}

In the sequel we call (DGCS)-phenomenon the instantaneous loss of
regularity described in point~(3) above.

\paragraph{\textmd{\emph{The dissipative case with constant 
coefficients}}}

If $\delta>0$ and $c(t)$ is a constant function (equal to 1 without loss of 
generality), equation (\ref{pbm:eqn}) reduces to
\begin{equation}
	u''(t)+2\delta A^{\sigma}u'(t)+Au(t)=0.
	\label{pbm:constant}
\end{equation}

Mathematical models with damping terms of this form were proposed
in~\cite{CR}, and then rigorously analyzed by many authors from
different points of view.  The first papers~\cite{CT1,CT2,CT3}, and
the more recent~\cite{HO}, are devoted to analyticity properties of
the semigroup associated to (\ref{pbm:constant}).  The classic
assumptions in these papers are that the operator $A$ is strictly
positive, $\sigma\in[0,1]$, and the phase space is $D(A^{1/2})\times
H$.  On a different side, the community working on dispersive
equations considered equation~(\ref{pbm:constant}) in the concrete
case where $\sigma\in[0,1]$ and $Au=-\Delta u$ in $\re^{n}$ or special
classes of unbounded domains.  They proved energy decay and dispersive
estimates, but exploiting in an essential way the spectral properties
of the Laplacian in those domains.  The interested reader is referred
to \cite{I1,I2,I3,shibata} and to the references quoted therein.

Finally, equation (\ref{pbm:constant}) was considered in~\cite{ggh:sd}
in full generality, namely for every $\sigma\geq 0$ and every
nonnegative self-adjoint operator $A$.  Two different regimes
appeared.  In the subcritical regime $\sigma\in[0,1/2]$, problem
(\ref{pbm:constant})--(\ref{pbm:data}) is well-posed in the classic
energy space $D(A^{1/2})\times H$ or more generally in
$D(A^{\alpha+1/2})\times D(A^{\alpha})$ with $\alpha\geq 0$.  In the
supercritical regime $\sigma\geq 1/2$, problem
(\ref{pbm:constant})--(\ref{pbm:data}) is well-posed in
$D(A^{\alpha})\times D(A^{\beta})$ if and only if
\begin{equation}
	1-\sigma\leq\alpha-\beta\leq\sigma.
	\label{hp:abs}
\end{equation}

This means that in the supercritical regime different choices of
the phase space are possible, even with $\alpha-\beta\neq 1/2$.

\paragraph{\textmd{\emph{The dissipative case with time-dependent 
coefficients}}}

As far as we know, the case of a dissipative equation with a
time-dependent propagation speed had not been considered yet.  The
main question we address in this paper is the extent to which the
dissipative term added in (\ref{pbm:eqn}) prevents the
(DGCS)-phenomenon of (\ref{pbm:non-diss}) from happening.  We discover
a composite picture, depending on $\sigma$.
\begin{itemize}
	\item In the subcritical regime $\sigma\in[0,1/2]$, if the strict
	hyperbolicity assumption~(\ref{hp:sh}) is satisfied,
	well-posedness results do depend on the time-regularity of $c(t)$
	(see Theorem~\ref{thm:sub-reg}).  Classic examples are the
	following.
	\begin{itemize}
		\item If $c(t)$ is $\alpha$-H\"{o}lder continuous for some
		exponent $\alpha>1-2\sigma$, then the dissipation prevails,
		and problem (\ref{pbm:eqn})--(\ref{pbm:data}) is well-posed in
		the classic energy space $D(A^{1/2})\times H$ or more
		generally in $D(A^{\beta+1/2})\times D(A^{\beta})$ with
		$\beta\geq 0$.
	
		\item If $c(t)$ is no more than $\alpha$-H\"{o}lder continuous
		for some exponent $\alpha<1-2\sigma$, then the dissipation can
		be neglected, so that (\ref{pbm:eqn}) behaves exactly as the
		non-dissipative equation (\ref{pbm:non-diss}).  This means
		well-posedness in the Gevrey space of order $(1-\alpha)^{-1}$
		and the possibility to produce the (DGCS)-phenomenon for less
		regular data (see Theorem~\ref{thm:dgcs}).
	
		\item The case with $\alpha=1-2\sigma$ is critical and also
		the size of the H\"{o}lder constant of $c(t)$ compared with
		$\delta$ comes into play.
	\end{itemize}

	\item In the supercritical regime $\sigma>1/2$ the dissipation
	prevails in an overwhelming way. In  Theorem~\ref{thm:sup-reg} we 
	prove that, if $c(t)$ is just measurable  and satisfies just the
	degenerate hyperbolicity condition
	\begin{equation}
		0\leq c(t)\leq\mu_{2},
		\label{hp:dh}
	\end{equation}
	then (\ref{pbm:eqn}) behaves as (\ref{pbm:constant}).  This means
	that problem (\ref{pbm:eqn})--(\ref{pbm:data}) is well-posed in
	$D(A^{\alpha})\times D(A^{\beta})$ if and only if (\ref{hp:abs})
	is satisfied, the same result obtained in~\cite{ggh:sd} in the
	case of a constant coefficient.

\end{itemize}

The second issue we address in this paper is the further
space-regularity of solutions for positive times, since a strong
dissipation is expected to have a regularizing effect similar to
parabolic equations.  This turns out to be true provided that the
assumptions of our well-posedness results are satisfied, and in
addition $\sigma\in(0,1)$.  Indeed, we prove that in this regime
$u(t)$ lies in the Gevrey space of order
$(2\min\{\sigma,1-\sigma\})^{-1}$ for every $t>0$.  We refer to
Theorem~\ref{thm:sup-gevrey} and Theorem~\ref{thm:sub-gevrey} for the
details.  This effect had already been observed in~\cite{MR:gevrey} in
the dispersive case.

We point out that the regularizing effect is maximum when $\sigma=1/2$
(the only case in which solutions become analytic with respect to
space variables) and disappears when $\sigma\geq 1$, meaning that a
stronger overdamping prevents smoothing. 

\paragraph{\textmd{\emph{Overview of the technique}}}

The spectral theory reduces the problem to an analysis of the family
of ordinary differential equations
\begin{equation}
	\ul''(t)+2\delta\lambda^{2\sigma}\ul'(t)+
	\lambda^{2}c(t)\ul(t)=0.
	\label{pbm:main-ode}
\end{equation}

When $\delta=0$, a coefficient $c(t)$ which oscillates with a suitable
period can produce a resonance effect so that (\ref{pbm:main-ode})
admits a solution whose oscillations have an amplitude which grows
exponentially with time.  This is the primordial origin of the
(DGCS)-phenomenon for non-dissipative equations.  When $\delta>0$, the
damping term causes an exponential decay of the amplitude of
oscillations.  \emph{The competition between the exponential energy
growth due to resonance and the exponential energy decay due to
dissipation originates the threshold effect we observed.}

When $c(t)$ is constant, equation (\ref{pbm:main-ode}) can be
explicitly integrated, and the explicit formulae for solutions led to
the sharp results of~\cite{ggh:sd}.  Here we need the same sharp
estimates, but without relying on explicit solutions.  To
this end, we introduce suitable energy estimates.  

In the supercritical regime $\sigma\geq 1/2$ we exploit the following
$\sigma$-adapted ``Kovaleskyan energy''
\begin{equation}
	E(t):=|\ul'(t)+\delta\lambda^{2\sigma}\ul(t)|^{2}+
	\delta^{2}\lambda^{4\sigma}|\ul(t)|^{2}.
	\label{energy:k}
\end{equation}

In the subcritical regime $\sigma\leq 1/2$ we exploit the so-called
``approximated hyperbolic energies''
\begin{equation}
	\Eep(t):=
	|\ul'(t)+\delta\lambda^{2\sigma}\ul(t)|^{2}+
	\delta^{2}\lambda^{4\sigma}|\ul(t)|^{2}+
	\lambda^{2}\cep(t)|\ul(t)|^{2},
	\label{energy:Eep}
\end{equation}
obtained by adding to (\ref{energy:k}) an ``hyperbolic term''
depending on a suitable smooth approximation $\cep(t)$ of $c(t)$,
which in turn is chosen in a $\lambda$-dependent way.  Terms of this
type are the key tool introduced in~\cite{dgcs} for the
non-dissipative equation.

\paragraph{\textmd{\emph{Future extensions}}}

We hope that this paper could represent a first step in the theory of
dissipative hyperbolic equations with variable coefficients, both
linear and nonlinear.  Next steps could be considering a coefficient
$c(x,t)$ depending both on time and space variables, and finally
quasilinear equations.  This could lead to improve the classic
results by K.~Nishihara~\cite{nishihara,nishihara-decay} for Kirchhoff
equations, whose linearization has a time-dependent coefficient, and
finally to consider more general local nonlinearities, in which case
the linearization involves a coefficient $c(x,t)$ depending on both
variables.

In a different direction, the subcritical case $\sigma\in[0,1/2]$ with
degenerate hyperbolicity assumptions remains open and could be the
subject of further research, in the same way as~\cite{cjs} was the
follow-up of~\cite{dgcs}.

On the other side, we hope that our counterexamples could finally
dispel the diffuse misconception according to which dissipative
hyperbolic equations are more stable, and hence definitely easier to
handle.  Now we know that a friction term below a suitable threshold
is substantially ineffective, opening the door to pathologies such as
the (DGCS)-phenomenon, exactly as in the non-dissipative case.

\paragraph{\textmd{\emph{Structure of the paper}}}
 
This paper is organized as follows.  In section~\ref{sec:notation} we
introduce the functional setting and we recall the classic existence
results from~\cite{dgcs}.  In section~\ref{sec:main} we state our main
results.  In section~\ref{sec:heuristics} we provide a heuristic
description of the competition between resonance and decay.  In
section~\ref{sec:proofs} we prove our existence and regularity
results.  In section~\ref{sec:counterexamples} we present our examples
of (DGCS)-phenomenon.

\setcounter{equation}{0}
\section{Notation and previous results}\label{sec:notation}

\paragraph{\textmd{\textit{Functional spaces}}}

Let $H$ be a separable Hilbert space.  Let us assume that $H$ admits a
countable complete orthonormal system $\{e_{k}\}_{k\in\n}$ made by
eigenvectors of $A$.  We denote the corresponding eigenvalues by
$\lambda_{k}^{2}$ (with the agreement that $\lk\geq 0$), so that
$Ae_{k}=\lambda_{k}^{2}e_{k}$ for every $k\in\n$.  In this case every
$u\in H$ can be written in a unique way in the form
$u=\sum_{k=0}^{\infty}u_{k}e_{k}$, where $u_{k}=\langle
u,e_{k}\rangle$ are the Fourier components of $u$.  In other words,
the Hilbert space $H$ can be identified with the set of sequences
$\{u_{k}\}$ of real numbers such that
$\sum_{k=0}^{\infty}u_{k}^{2}<+\infty$.

We stress that this is just a simplifying assumption, with
substantially no loss of generality.  Indeed, according to the
spectral theorem in its general form (see for example Theorem~VIII.4
in~\cite{reed}), one can always identify $H$ with $L^{2}(M,\mu)$ for a
suitable measure space $(M,\mu)$, in such a way that under this
identification the operator $A$ acts as a multiplication operator by
some measurable function $\lambda^{2}(\xi)$.  All definitions and
statements in the sequel, with the exception of the counterexamples
of Theorem~\ref{thm:dgcs}, can be easily extended to the
general setting just by replacing the sequence $\{\lk^{2}\}$ with the
function $\lambda^{2}(\xi)$, and the sequence $\{u_{k}\}$ of Fourier
components of $u$ with the element $\widehat{u}(\xi)$ of
$L^{2}(M,\mu)$ corresponding to $u$ under the identification of $H$
with $L^{2}(M,\mu)$.

The usual functional spaces can be characterized in terms of Fourier
components as follows.

\begin{defn}
	\begin{em}
		Let $u$ be a sequence $\{u_{k}\}$ of real numbers.
		\begin{itemize}
			\item  \emph{Sobolev spaces}. For every $\alpha\geq 0$ it 
			turns out that $u\in D(A^{\alpha})$ if 
			\begin{equation}
				\|u\|_{D(A^{\alpha})}^{2}:=
				\sum_{k=0}^{\infty}(1+\lk)^{4\alpha}u_{k}^{2}<+\infty.
				\label{defn:sobolev}
			\end{equation}
		
			\item \emph{Distributions}.  We say that $u\in
			D(A^{-\alpha})$ for some $\alpha\geq 0$ if 
			\begin{equation}
				\|u\|_{D(A^{-\alpha})}^{2}:=
				\sum_{k=0}^{\infty}(1+\lk)^{-4\alpha}u_{k}^{2}<+\infty.
				\label{defn:distributions}
			\end{equation}
		
			\item \emph{Generalized Gevrey spaces}.  Let
			$\varphi:[0,+\infty)\to[0,+\infty)$ be any function, let
			$r\geq 0$, and let $\alpha\in\re$.  We say that
			$u\in\G_{\varphi,r,\alpha}(A)$ if 
			\begin{equation}
				\|u\|^{2}_{\varphi,r,\alpha}:=
				\sum_{k=0}^{\infty}(1+\lk)^{4\alpha}
				u_{k}^{2} \exp\left(\strut 2r\varphi(\lk)\right)
				<+\infty.
				\label{defn:gevrey}
			\end{equation}
		
			\item \emph{Generalized Gevrey ultradistributions}.  Let
			$\psi:[0,+\infty)\to[0,+\infty)$ be any function, let
			$R\geq 0$, and let $\alpha\in\re$.  We say that
			$u\in\G_{-\psi,R,\alpha}(A)$ if 
			\begin{equation}
				\|u\|^{2}_{-\psi,R,\alpha}:=
				\sum_{k=0}^{\infty}(1+\lk)^{4\alpha}
				u_{k}^{2} \exp\left(\strut -2R\psi(\lk)\right)
				<+\infty.
				\label{defn:hyperfunctions}
			\end{equation}

		\end{itemize}
	\end{em}
\end{defn}

\begin{rmk}
	\begin{em}
		If $\varphi_{1}(x)=\varphi_{2}(x)$ for every $x>0$, then
		$\G_{\varphi_{1},r,\alpha}(A)=\G_{\varphi_{2},r,\alpha}(A)$
		for every admissible value of $r$ and $\alpha$.  For this
		reason, with a little abuse of notation, we consider the
		spaces $\G_{\varphi,r,\alpha}(A)$ even when $\varphi(x)$ is
		defined only for $x>0$.  The same comment applies also to the
		spaces $\G_{-\psi,R,\alpha}(A)$.
	\end{em}
\end{rmk}

The quantities defined in (\ref{defn:sobolev}) through
(\ref{defn:hyperfunctions}) are actually norms which induce a Hilbert
space structure on $D(A^{\alpha})$, $\G_{\varphi,r,\alpha}(A)$,
$\G_{-\psi,R,\alpha}(A)$, respectively.  The standard inclusions
$$\G_{\varphi,r,\alpha}(A)\subseteq D(A^{\alpha})\subseteq H
\subseteq D(A^{-\alpha})\subseteq\G_{-\psi,R,-\alpha}(A)$$
hold true for every $\alpha\geq 0$ and every admissible choice of
$\varphi$, $\psi$, $r$, $R$.  All inclusions are strict if $\alpha$,
$r$ and $R$ are positive, and the sequences $\{\lk\}$,
$\{\varphi(\lk)\}$, and $\{\psi(\lk)\}$ are unbounded.

We observe that $\G_{\varphi,r,\alpha}(A)$ is actually a so-called
\emph{scale of Hilbert spaces} with respect to the parameter $r$, with
larger values of $r$ corresponding to smaller spaces.  Analogously,
$\G_{-\psi,R,\alpha}(A)$ is a scale of Hilbert spaces with
respect to the parameter $R$, but with larger values of $R$
corresponding to larger spaces.

\begin{rmk}
	\begin{em}
		Let us consider the concrete case where $I\subseteq\re$ is an
		open interval, $H=L^{2}(I)$, and $Au=-u_{xx}$, with periodic
		boundary conditions.  For every $\alpha\geq 0$, the space
		$D(A^{\alpha})$ is actually the usual Sobolev space
		$H^{2\alpha}(I)$, and $D(A^{-\alpha})$ is the usual space of
		distributions of order $2\alpha$.
		
		When $\varphi(x):=x^{1/s}$ for some $s>0$, elements of
		$\G_{\varphi,r,0}(A)$ with $r>0$ are usually called Gevrey
		functions of order $s$, the case $s=1$ corresponding to
		analytic functions.  When $\psi(x):=x^{1/s}$ for some $s>0$,
		elements of $\G_{-\psi,R,0}(A)$ with $R>0$ are usually called
		Gevrey ultradistributions of order $s$, the case $s=1$
		corresponding to analytic functionals.  In this case the
		parameter $\alpha$ is substantially irrelevant because the
		exponential term is dominant both in (\ref{defn:gevrey}) and
		in (\ref{defn:hyperfunctions}).
		
		For the sake of consistency, with a little abuse of notation
		we use the same terms (Gevrey functions, Gevrey
		ultradistributions, analytic functions and analytic
		functionals) in order to denote the same spaces also in the
		general abstract framework.  To be more precise, we should
		always add ``with respect to the operator $A$'', or even
		better ``with respect to the operator $A^{1/2}$''.
	\end{em}
\end{rmk}

\paragraph{\textmd{\textit{Continuity moduli}}}

Throughout this paper we call \emph{continuity modulus} any continuous
function $\omega:[0,+\infty)\to[0,+\infty)$ such that $\omega(0)=0$, 
$\omega(x)>0$ for every $x>0$, and moreover
\begin{equation}
	x\to\omega(x)\mbox{ is a nondecreasing function},
	\label{hp:omega-monot-0}
\end{equation}
\begin{equation}
	x\to\frac{x}{\omega(x)}\mbox{ is a nondecreasing function.}
	\label{hp:omega-monot}
\end{equation}

A function $c:[0,+\infty)\to\re$ is said to be $\omega$-continuous if
\begin{equation}
	|c(a)-c(b)|\leq
	\omega(|a-b|)
	\quad\quad
	\forall a\geq 0,\ \forall b\geq 0.
	\label{hp:ocont}
\end{equation}

More generally, a function $c:X\to\re$ (with $X\subseteq\re$) is said
to be $\omega$-continuous if it satisfies the same inequality for
every $a$ and $b$ in $X$.

\paragraph{\textmd{\textit{Previous results}}}\label{sec:thmbibl}

We are now ready to recall the classic results concerning existence,
uniqueness, and regularity for solutions to problem
(\ref{pbm:eqn})--(\ref{pbm:data}).  We state them using our notations
which allow general continuity moduli and general spaces of Gevrey
functions or ultradistributions.

Proofs are a straightforward application of the approximated energy
estimates introduced in~\cite{dgcs}.  In that paper only the case
$\delta=0$ is considered, but when $\delta\geq 0$ all new terms have
the ``right sign'' in those estimates.

The first result concerns existence and uniqueness in huge spaces
such as analytic functionals, with minimal assumptions on $c(t)$.

\begin{thmbibl}[see~{\cite[Theorem~1]{dgcs}}]\label{thmbibl:dh}
	Let us consider problem (\ref{pbm:eqn})--(\ref{pbm:data}) under 
	the following assumptions:
	\begin{itemize}
		\item $A$ is a self-adjoint nonnegative operator on a
		separable Hilbert space $H$,
	
		\item $c\in L^{1}((0,T))$ for every $T>0$ (without sign 
		conditions),
	
		\item $\sigma\geq 0$ and $\delta\geq 0$ are two real numbers,
	
		\item initial conditions satisfy 
		$$(u_{0},u_{1})\in
		\G_{-\psi,R_{0},1/2}(A)\times\G_{-\psi,R_{0},0}(A)$$ 
		for some
		$R_{0}>0$ and some $\psi:(0,+\infty)\to(0,+\infty)$ such that
		$$\limsup_{x\to +\infty}\frac{x}{\psi(x)}<+\infty.$$

	\end{itemize}
	
	Then there exists a nondecreasing function
	$R:[0,+\infty)\to[0,+\infty)$, with $R(0)=R_{0}$, such that problem
	(\ref{pbm:eqn})--(\ref{pbm:data}) admits a unique solution
	\begin{equation}
		u\in
		C^{0}\left([0,+\infty);\G_{-\psi,R(t),1/2}(A)\right)
		\cap
		C^{1}\left([0,+\infty);\G_{-\psi,R(t),0}(A)\right).		
		\label{th:(DGCS)-psi}
	\end{equation}
\end{thmbibl}

Condition (\ref{th:(DGCS)-psi}), with the range space increasing with time,
simply means that 
$$u\in
C^{0}\left([0,\tau];\G_{-\psi,R(\tau),1/2}(A)\right)\cap
C^{1}\left([0,\tau];\G_{-\psi,R(\tau),0}(A)\right)
\quad\quad
\forall\tau\geq 0.$$
		
This amounts to say that scales of Hilbert spaces, rather than fixed
Hilbert spaces, are the natural setting for this problem.

In the second result we assume strict hyperbolicity and 
$\omega$-continuity of the coefficient, and we obtain well-posedness 
in a suitable class of Gevrey ultradistributions.

\begin{thmbibl}[see~{\cite[Theorem~3]{dgcs}}]\label{thmbibl:existence}
	Let us consider problem (\ref{pbm:eqn})--(\ref{pbm:data}) under 
	the following assumptions:
	\begin{itemize}
		\item $A$ is a self-adjoint nonnegative operator on a
		separable Hilbert space $H$,
	
		\item  the coefficient $c:[0,+\infty)\to\re$ satisfies the 
		strict hyperbolicity assumption (\ref{hp:sh}) and the 
		$\omega$-continuity assumption (\ref{hp:ocont}) for some 
		continuity modulus $\omega(x)$,
	
		\item $\sigma\geq 0$ and $\delta\geq 0$ are two real numbers,
	
		\item initial conditions satisfy
		$$(u_{0},u_{1})\in
		\G_{-\psi,R_{0},1/2}(A)\times\G_{-\psi,R_{0},0}(A)$$ 
		for some $R_{0}>0$ and some function
		$\psi:(0,+\infty)\to(0,+\infty)$ such that
		\begin{equation}
			\limsup_{x\to +\infty}\frac{x}{\psi(x)}\,
			\omega\left(\frac{1}{x}\right)<+\infty.
			\label{hp:(DGCS)-psi}
		\end{equation}

	\end{itemize}
	
	Let $u$ be the unique solution to the problem provided by 
	Theorem~\ref{thmbibl:dh}. 
	
	Then there exists $R>0$ such that 
	$$u\in C^{0}\left([0,+\infty),\G_{-\psi,R_{0}+Rt,1/2}(A)\right)
	\cap C^{1}\left([0,+\infty),\G_{-\psi,R_{0}+Rt,0}(A)\right).$$
	
\end{thmbibl}

The third result we recall concerns existence of regular solutions.
The assumptions on $c(t)$ are the same as in
Theorem~\ref{thmbibl:existence}, but initial data are significantly
more regular (Gevrey spaces instead of Gevrey ultradistributions).

\begin{thmbibl}[see~{\cite[Theorem~2]{dgcs}}]\label{thmbibl:regularity}
	Let us consider problem (\ref{pbm:eqn})--(\ref{pbm:data}) under 
	the following assumptions:
	\begin{itemize}
		\item $A$ is a self-adjoint nonnegative operator on a
		separable Hilbert space $H$,
	
		\item  the coefficient $c:[0,+\infty)\to\re$ satisfies the 
		strict hyperbolicity assumption (\ref{hp:sh}) and the 
		$\omega$-continuity assumption (\ref{hp:ocont}) for some 
		continuity modulus $\omega(x)$,
	
		\item $\sigma\geq 0$ and $\delta\geq 0$ are two real numbers,
	
		\item initial conditions satisfy
		$$(u_{0},u_{1})\in
		\G_{\varphi,r_{0},1/2}(A)\times\G_{\varphi,r_{0},0}(A)$$ 
		for some $r_{0}>0$ and some function
		$\varphi:(0,+\infty)\to(0,+\infty)$ such that
		\begin{equation}
			\limsup_{x\to +\infty}\frac{x}{\varphi(x)}\,
			\omega\left(\frac{1}{x}\right)<+\infty.
			\label{hp:(DGCS)-phi}
		\end{equation}

	\end{itemize}
	
	Let $u$ be the unique solution to the problem provided by 
	Theorem~\ref{thmbibl:dh}. 
	
	Then there exist $T>0$ and $r>0$ such that $rT<r_{0}$ and
	\begin{equation}
		u\in
		C^{0}\left([0,T],\G_{\varphi,r_{0}-rt,1/2}(A)\right)\cap
		C^{1}\left([0,T],\G_{\varphi,r_{0}-rt,0}(A)\right).
		\label{th:(DGCS)-phi}
	\end{equation}
	
\end{thmbibl}

\begin{rmk}
	\begin{em}
		The key assumptions of Theorem~\ref{thmbibl:existence} and
		Theorem~\ref{thmbibl:regularity} are (\ref{hp:(DGCS)-psi}) and
		(\ref{hp:(DGCS)-phi}), respectively, representing the exact
		compensation between space-regularity of initial data and
		time-regularity of the coefficient $c(t)$ required in order to
		obtain well-posedness.  
		
		These conditions do not appear explicitly in~\cite{dgcs},
		where they are replaced by suitable specific choices of
		$\omega$, $\varphi$, $\psi$, which of course satisfy the same
		relations.  To our knowledge, those conditions were stated for
		the first time in~\cite{gg:der-loss}, thus unifying several
		papers that in the last 30 years had been devoted to special
		cases (see for example~\cite{colombini} and the references
		quoted therein).
	\end{em}
\end{rmk}

\begin{rmk}\label{rmk:dgcs-1}
	\begin{em}
		The standard example of application of
		Theorem~\ref{thmbibl:existence} and
		Theorem~\ref{thmbibl:regularity} is the following.  Let us
		assume that $c(t)$ is $\alpha$-H\"{o}lder continuous for some
		$\alpha\in(0,1)$, namely $\omega(x)=Mx^{\alpha}$ for a
		suitable constant $M$.  Then (\ref{hp:(DGCS)-psi}) and
		(\ref{hp:(DGCS)-phi}) hold true with
		$\psi(x)=\varphi(x):=x^{1-\alpha}$.  This leads to
		well-posedness both in the large space of Gevrey
		ultradistributions of order $(1-\alpha)^{-1}$, and in the
		small space of Gevrey functions of the same order.
	\end{em}
\end{rmk}

\begin{rmk}\label{rmk:dgcs-2}
	\begin{em}
		The choice of ultradistributions in
		Theorem~\ref{thmbibl:existence} is not motivated by the search
		for generality, but it is in some sense the only possible one
		because of the (DGCS)-phenomenon exhibited in~\cite{dgcs}, at
		least in the non-dissipative case.  When $\delta=0$, if
		initial data are taken in Sobolev spaces or in any space
		larger than the Gevrey spaces of
		Theorem~\ref{thmbibl:regularity}, then it may happen that for
		all positive times the solution lies in the space of
		ultradistributions specified in
		Theorem~\ref{thmbibl:existence}, and nothing more.  In other
		words, for $\delta=0$ there is no well-posedness result in
		between the Gevrey spaces of Theorem~\ref{thmbibl:regularity}
		and the Gevrey ultradistributions of
		Theorem~\ref{thmbibl:existence}, and conditions
		(\ref{hp:(DGCS)-psi}) and (\ref{hp:(DGCS)-phi}) are optimal.
		
		The aim of this paper is to provide an optimal picture for the
		case $\delta>0$.
	\end{em}
\end{rmk}

\setcounter{equation}{0}
\section{Main results}\label{sec:main}

In this section we state our main regularity results for solutions to 
(\ref{pbm:eqn})--(\ref{pbm:data}). To this end, we need some further 
notation. Given any $\nu\geq 0$, we write $H$ as an orthogonal direct sum
\begin{equation}
	H:=H_{\nu,-}\oplus H_{\nu,+},
	\label{defn:H+-}
\end{equation}
where $H_{\nu,-}$ is the closure of the subspace generated by all
eigenvectors of $A$ relative to eigenvalues $\lk<\nu$, and $H_{\nu,+}$
is the closure of the subspace generated by all eigenvectors of $A$
relative to eigenvalues $\lk\geq\nu$.  For every vector $u\in H$, we
write $u_{\nu,-}$ and $u_{\nu,+}$ to denote its components with
respect to the decomposition (\ref{defn:H+-}).  We point out that
$H_{\nu,-}$ and $H_{\nu,+}$ are $A$-invariant subspaces of $H$, and
that $A$ is a bounded operator when restricted to $H_{\nu,-}$, and a
coercive operator when restricted to $H_{\nu,+}$ if $\nu>0$.

In the following statements we provide separate estimates for
low-frequency components $u_{\nu,-}(t)$ and high-frequency
components $u_{\nu,+}(t)$ of solutions to (\ref{pbm:eqn}).  This is
due to the fact that the energy of $u_{\nu,-}(t)$ can be unbounded as
$t\to +\infty$, while in many cases we are able to prove that the
energy of $u_{\nu,+}(t)$ is bounded in time.

\subsection{Existence results in Sobolev spaces}\label{sec:main-existence}

The first result concerns the supercritical regime
$\sigma\geq 1/2$, in which case the dissipation always dominates the 
time-dependent coefficient.

\begin{thm}[Supercritical dissipation]\label{thm:sup-reg}
	Let us consider problem (\ref{pbm:eqn})--(\ref{pbm:data}) under 
	the following assumptions:
	\begin{itemize}
		\item $A$ is a self-adjoint nonnegative operator on a
		separable Hilbert space $H$,
	
		\item the coefficient $c:[0,+\infty)\to\re$ is measurable and
		satisfies the degenerate hyperbolicity assumption
		(\ref{hp:dh}),
	
		\item $\sigma$ and $\delta$ are two positive real numbers such
		that either $\sigma>1/2$, or $\sigma=1/2$ and
		$4\delta^{2}\geq\mu_{2}$,
		
		\item $(u_{0},u_{1})\in D(A^{\alpha})\times D(A^{\beta})$
		for some real numbers $\alpha$ and $\beta$ satisfying (\ref{hp:abs}).
	\end{itemize}

	Let $u$ be the unique solution to the problem provided by 
	Theorem~\ref{thmbibl:dh}. 
	
	Then $u$ actually satisfies
	\begin{equation}
		(u,u')\in C^{0}\left([0,+\infty),D(A^{\alpha})\times D(A^{\beta})\right).
		\label{th:sup-reg}
	\end{equation}
	
	Moreover, for every $\nu\geq 1$ such that
	$4\delta^{2}\nu^{4\sigma-2}\geq\mu_{2}$, it turns out that
	\begin{equation}
		|A^{\beta}u_{\nu,+}'(t)|^{2}+|A^{\alpha}u_{\nu,+}(t)|^{2}\leq 
		\left(2+\frac{2}{\delta^{2}}+\frac{\mu_{2}^{2}}{\delta^{4}}\right)
		|A^{\beta}u_{1,\nu,+}|^{2}+
		3\left(1+\frac{\mu_{2}^{2}}{2\delta^{2}}\right)
		|A^{\alpha}u_{0,\nu,+}|^{2}
		\label{th:sup-u+}
	\end{equation}
	for every $t\geq 0$.
\end{thm}

Our second result concerns the subcritical regime $\sigma\in[0,1/2]$, 
in which case the time-regularity of $c(t)$ competes with the 
exponent $\sigma$.

\begin{thm}[Subcritical dissipation]\label{thm:sub-reg}
	Let us consider problem (\ref{pbm:eqn})--(\ref{pbm:data}) under 
	the following assumptions:
	\begin{itemize}
		\item $A$ is a self-adjoint nonnegative operator on a
		separable Hilbert space $H$,
	
		\item  the coefficient $c:[0,+\infty)\to\re$ satisfies the 
		strict hyperbolicity assumption (\ref{hp:sh}) and the 
		$\omega$-continuity assumption (\ref{hp:ocont}) for some 
		continuity modulus $\omega(x)$,
	
		\item $\sigma\in[0,1/2]$ and $\delta>0$ are two real 
		numbers such that
		\begin{equation}
			4\delta^{2}\mu_{1}>\Lambda_{\infty}^{2}+2\delta\Lambda_{\infty},
			\label{hp:sub}
		\end{equation}
		where we set
		\begin{equation}
			\Lambda_{\infty}:=
			\limsup_{\ep\to 0^{+}}\frac{\omega(\ep)}{\ep^{1-2\sigma}},
			\label{defn:Linfty}
		\end{equation}
	
		\item $(u_{0},u_{1})\in D(A^{1/2})\times H$.

	\end{itemize}
	
	Let $u$ be the unique solution to the problem provided by 
	Theorem~\ref{thmbibl:dh}. 
	
	Then $u$ actually satisfies
	$$u\in C^{0}\left([0,+\infty),D(A^{1/2})\right) \cap
	C^{1}\left([0,+\infty),H\right).$$
		
	Moreover, for every $\nu\geq 1$ such that
	\begin{equation}
		4\delta^{2}\mu_{1}\geq
		\left[\lambda^{1-2\sigma}\omega\left(\frac{1}{\lambda}\right)\right]^{2}
		+2\delta
		\left[\lambda^{1-2\sigma}\omega\left(\frac{1}{\lambda}\right)\right]
		\label{hp:sub-lambda}
	\end{equation}
	for every $\lambda\geq\nu$, it turns out that
	\begin{equation}
		|u_{\nu,+}'(t)|^{2}+2\mu_{1}|A^{1/2}u_{\nu,+}(t)|^{2}\leq
		4|u_{1,\nu,+}|^{2}+2(3\delta^{2}+
		\mu_{2})|A^{1/2}u_{0,\nu,+}|^{2}
		\label{th:sub-u+}
	\end{equation}
	for every $t\geq 0$.
\end{thm}

Let us make a few comments on the first two statements.

\begin{rmk}\label{rmk:low-freq}
	\begin{em}
		In both results we proved that a suitable high-frequency
		component of the solution can be uniformly bounded in
		terms of initial data.  Low-frequency components might in
		general diverge as $t\to +\infty$.  Nevertheless, they can
		always be estimated as follows.
		
		Let us just assume that $c\in L^{1}((0,T))$ for every $T>0$.
		Then for every $\nu\geq 0$ the component $u_{\nu,-}(t)$
		satisfies
		\begin{equation}
			|u_{\nu,-}'(t)|^{2}+|A^{1/2}u_{\nu,-}(t)|^{2}\leq 
			\left(|u_{1,\nu,-}|^{2}+|A^{1/2}u_{0,\nu,-}|^{2}\right)
			\exp\left(\nu t+\nu\int_{0}^{t}|c(s)|\,ds\right)
			\label{th:est-u-}
		\end{equation}
		for every $t\geq 0$.  Indeed, let $F(t)$ denote the left-hand
		side of (\ref{th:est-u-}). Then
		\begin{eqnarray*}
			F'(t) & = & -4\delta|A^{\sigma/2}u_{\nu,-}'(t)|^{2}+
			2(1-c(t))\langle u_{\nu,-}'(t),Au_{\nu,-}(t)\rangle \\
			 & \leq & 2(1+|c(t)|)\cdot
			 |u_{\nu,-}'(t)|\cdot\nu|A^{1/2}u_{\nu,-}(t)|  \\
			 & \leq & \nu(1+|c(t)|) F(t)
		\end{eqnarray*}
		for almost every $t\geq 0$, so that (\ref{th:est-u-}) follows
		by integrating this differential inequality.
	\end{em}
\end{rmk}

\begin{rmk}
	\begin{em}
		The phase spaces involved in Theorem~\ref{thm:sup-reg} and
		Theorem~\ref{thm:sub-reg} are exactly the same which are known
		to be optimal when $c(t)$ is constant (see~\cite{ggh:sd}).  In
		particular, the only possible choice in the subcritical regime
		is the classic energy space $D(A^{1/2})\times H$, or more
		generally $D(A^{\alpha+1/2})\times D(A^{\alpha})$.  This ``gap
		1/2'' between the powers of $A$ involved in the phase space is
		typical of hyperbolic problems, and it is the same which
		appears in the classic results of
		section~\ref{sec:notation}.
		
		On the contrary, in the supercritical regime there is  
		an interval of possible gaps, described by~(\ref{hp:abs}). This
		interval is always centered in 1/2, but also different values
		are allowed, including negative ones when $\sigma>1$.
	\end{em}
\end{rmk}

\begin{rmk}
	\begin{em}
		The classic example of application of
		Theorem~\ref{thm:sub-reg} is the following.  Let us assume
		that $c(t)$ is $\alpha$-H\"{o}lder continuous for some
		$\alpha\in(0,1)$, namely $\omega(x)=Mx^{\alpha}$ for some
		constant $M$.  Then problem (\ref{pbm:eqn})--(\ref{pbm:data})
		is well-posed in the energy space provided that either
		$\alpha>1-2\sigma$, or $\alpha=1-2\sigma$ and $M$ is small
		enough.  Indeed, for $\alpha>1-2\sigma$ we get
		$\Lambda_{\infty}=0$, and hence (\ref{hp:sub}) is
		automatically satisfied.  For $\alpha=1-2\sigma$ we get
		$\Lambda_{\infty}=M$, so that (\ref{hp:sub}) is satisfied
		provided that $M$ is small enough.
		
		In all other cases, namely when either $\alpha<1-2\sigma$, or
		$\alpha=1-2\sigma$ and $M$ is large enough, only
		Theorem~\ref{thmbibl:existence} applies to initial data in
		Sobolev spaces, providing global existence just in the sense
		of Gevrey ultradistributions of order $(1-\alpha)^{-1}$.
	\end{em}
\end{rmk}

\begin{rmk}\label{rmk:lip-0}
	\begin{em}
		Let us examine the limit case $\sigma=0$, which falls in the
		subcritical regime.  
		
		When $\sigma=0$, assumption (\ref{hp:sub}) is satisfied if and
		only if $c(t)$ is Lipschitz continuous and its Lipschitz
		constant is small enough.  On the other hand, in the Lipschitz
		case it is a classic result that problem
		(\ref{pbm:eqn})--(\ref{pbm:data}) is well-posed in the energy
		space, regardless of the Lipschitz constant.  Therefore, the
		result stated in Theorem~\ref{thm:sub-reg} is non-optimal when
		$\sigma=0$ and $c(t)$ is Lipschitz continuous.
		
		A simple refinement of our argument would lead to the full
		result also in this case, but unfortunately it would be
		useless in all other limit cases in which $c(t)$ is
		$\alpha$-H\"{o}lder continuous with $\alpha=1-2\sigma$ and
		$\sigma\in(0,1/2]$.  We refer to section~\ref{sec:heuristics}
		for further details.
	\end{em}
\end{rmk}

\begin{rmk}\label{rmk:1/2}
	\begin{em}
		Let us examine the limit case $\sigma=1/2$, which falls both
		in the subcritical and in the supercritical regime, so that
		the conclusions of Theorem~\ref{thm:sup-reg} and
		Theorem~\ref{thm:sub-reg} coexist. Both of them provide 
		well-posedness in the energy space, but with different 
		assumptions.
		
		Theorem~\ref{thm:sup-reg} needs less assumptions on $c(t)$,
		which is only required to be measurable and to satisfy the
		degenerate hyperbolicity assumption~(\ref{hp:dh}), but it
		requires $\delta$ to be large enough so that
		$4\delta^{2}\geq\mu_{2}$.
		
		On the contrary, Theorem~\ref{thm:sub-reg} needs less
		assumptions on $\delta$, which is only required to be
		positive, but it requires $c(t)$ to be continuous and to
		satisfy the strict hyperbolicity assumption~(\ref{hp:sh}).
		Indeed, inequality (\ref{hp:sub}) is automatically satisfied
		in the case $\sigma=1/2$ because $\Lambda_{\infty}=0$.
		
		The existence of two different sets of assumptions leading to
		the same conclusion suggests the existence of a unifying
		statement, which could probably deserve further investigation.
		
	\end{em}
\end{rmk}

\subsection{Gevrey regularity for positive times}\label{sec:gevrey}

A strong dissipation in the range $\sigma\in(0,1)$ has a
regularizing effect on initial data, provided that the solution 
exists in Sobolev spaces.  In the following two statements
we quantify this effect in terms of scales of Gevrey spaces.

Both results can be summed up by saying that the solution lies, for 
positive times, in Gevrey spaces of order 
$(2\min\{\sigma,1-\sigma\})^{-1}$. It is not difficult to show that 
this order is optimal, even in the case where $c(t)$ is constant.

\begin{thm}[Supercritical dissipation]\label{thm:sup-gevrey}
	Let us consider problem (\ref{pbm:eqn})--(\ref{pbm:data}) under 
	the same assumptions of Theorem~\ref{thm:sup-reg}, and let $u$ be
	the unique solution to the problem provided by
	Theorem~\ref{thmbibl:dh}.
	
	Let us assume in addition that either $\sigma\in(1/2,1)$, or
	$\sigma=1/2$ and $4\delta^{2}>\mu_{2}$.  Let us set
	$\varphi(x):=x^{2(1-\sigma)}$, and
	\begin{equation}
		C(t):=\int_{0}^{t}c(s)\,ds.
		\label{defn:C(t)}
	\end{equation}
	
	Then there exists $r>0$ such that
	\begin{equation}
		(u,u')\in C^{0}\left((0,+\infty),\mathcal{G}_{\varphi,\alpha,rC(t)}(A)
		\times\mathcal{G}_{\varphi,\beta,rC(t)}(A)\right),
		\label{th:sup-gevrey}
	\end{equation}
	and there exist $\nu\geq 1$ and $K>0$ such that
	\begin{equation}
		\|u_{\nu,+}'(t)\|_{\varphi,\beta,rC(t)}^{2}+
		\|u_{\nu,+}(t)\|_{\varphi,\alpha,rC(t)}^{2}
		\leq K\left(
		|A^{\beta}u_{1,\nu,+}|^{2}+|A^{\alpha}u_{0,\nu,+}|^{2}
		\right)
		\label{th:sup-g-est}
	\end{equation}
	for every $t>0$.  The constants $r$, $\nu$, and $K$ depend only on
	$\delta$, $\mu_{2}$, and $\sigma$.
	
\end{thm}

Of course, (\ref{th:sup-gevrey}) and (\ref{th:sup-g-est}) are
nontrivial only if $C(t)>0$, which is equivalent to saying that the
coefficient $c(t)$ is not identically 0 in $[0,t]$.  On the other
hand, this weak form of hyperbolicity is necessary, since no
regularizing effect on $u(t)$ can be expected as long as $c(t)$
vanishes.

\begin{thm}[Subcritical dissipation]\label{thm:sub-gevrey}
	Let us consider problem (\ref{pbm:eqn})--(\ref{pbm:data}) under
	the same assumptions of Theorem~\ref{thm:sub-reg}, and let $u$ be
	the unique solution to the problem provided by
	Theorem~\ref{thmbibl:dh}.
	
	Let us assume in addition that $\sigma\in(0,1/2]$ (instead of
	$\sigma\in[0,1/2]$), and let us set $\varphi(x):=x^{2\sigma}$.
	
	Then there exists $r>0$ such that
	$$u\in
	C^{0}\left((0,+\infty),\mathcal{G}_{\varphi,1/2,rt}(A)\right) \cap
	C^{1}\left((0,+\infty),\mathcal{G}_{\varphi,0,rt}(A)\right),$$
	and there exist $\nu\geq 1$ and $K>0$ such that
	\begin{equation}
		\|u_{\nu,+}'(t)\|_{\varphi,0,rt}^{2}+
		\|u_{\nu,+}(t)\|_{\varphi,1/2,rt}^{2}
		\leq K
		\left(|u_{1,\nu,+}|^{2}+|A^{1/2}u_{0,\nu,+}|^{2}\right)
		\label{th:sub-g-est}
	\end{equation}
	for every $t>0$.  The constants $r$, $\nu$, and $K$ depend only on
	$\delta$, $\mu_{1}$, $\mu_{2}$, $\sigma$ and $\omega$.
\end{thm}

The estimates which provide Gevrey regularity of high-frequency 
components provide also the decay of the same components as $t\to 
+\infty$. We refer to Lemma~\ref{lemma:ODE-sup} and 
Lemma~\ref{lemma:ODE-sub} for further details.

\subsection{Counterexamples}

The following result shows that even strongly dissipative hyperbolic 
equations can exhibit the (DGCS)-phenomenon, provided that we are 
in the subcritical regime.

\begin{thm}[(DGCS)-phenomenon]\label{thm:dgcs}
	Let $A$ be a linear operator on a Hilbert space $H$. 
	Let us assume that there exists a countable (not necessarily 
	complete) orthonormal system $\{e_{k}\}$ in $H$, and an 
	unbounded sequence $\{\lk\}$ of positive real numbers 
	such that $Ae_{k}=\lk^{2}e_{k}$ for every $k\in\n$. Let 
	$\sigma\in[0,1/2)$ and $\delta>0$ be real numbers.
	
	Let $\omega:[0,+\infty)\to[0,+\infty)$ be a continuity modulus 
	such that
	\begin{equation}
		\lim_{\ep\to 0^{+}}\frac{\omega(\ep)}{\ep^{1-2\sigma}}=+\infty.
		\label{hp:cex-omega}
	\end{equation}
	
	Let $\varphi:(0,+\infty)\to(0,+\infty)$ and 
	$\psi:(0,+\infty)\to(0,+\infty)$ be two functions such that
	\begin{equation}
		\lim_{x\to+\infty}\frac{x}{\varphi(x)}\,\omega\left(\frac{1}{x}\right)=
		\lim_{x\to+\infty}\frac{x}{\psi(x)}\,\omega\left(\frac{1}{x}\right)=
		+\infty.
		\label{hp:cex-phi-psi}
	\end{equation}
	
	Then there exist a function $c:\re\to\re$ such that
	\begin{equation}
		\frac{1}{2}\leq c(t)\leq\frac{3}{2}
		\quad\quad
		\forall t\in\re,
		\label{th:c-bound}
	\end{equation}
	\begin{equation}
		|c(t)-c(s)|\leq\omega(|t-s|)
		\quad\quad
		\forall(t,s)\in\re^{2},
		\label{th:c-omega}
	\end{equation}
	and a solution $u(t)$ to equation (\ref{pbm:eqn}) such that
	\begin{equation}
		(u(0),u'(0))\in
		\mathcal{G}_{\varphi,r,1/2}(A)\times\mathcal{G}_{\varphi,r,0}(A)
		\quad\quad
		\forall r>0,
		\label{th:u0}
	\end{equation}
	\begin{equation}
		(u(t),u'(t))\not\in
		\mathcal{G}_{-\psi,R,1/2}(A)\times\mathcal{G}_{-\psi,R,0}(A)
		\quad\quad
		\forall R>0,\ \forall t>0.
		\label{th:ut}
	\end{equation}
\end{thm}

\begin{rmk}
	\begin{em}
		Due to (\ref{th:c-bound}), (\ref{th:c-omega}), and
		(\ref{th:u0}), the function $u(t)$ provided by
		Theorem~\ref{thm:dgcs} is a solution to (\ref{pbm:eqn}) in the
		sense of Theorem~\ref{thmbibl:dh} with $\psi(x):=x$, or
		even better in the sense of Theorem~\ref{thmbibl:existence}
		with $\psi(x):=x\omega(1/x)$.
	\end{em}
\end{rmk}

\begin{rmk}\label{rmk:L-infty}
	\begin{em}
		Assumption (\ref{hp:cex-omega}) is equivalent to saying that
		$\Lambda_{\infty}$ defined by (\ref{defn:Linfty}) is equal to
		$+\infty$, so that (\ref{hp:sub}) can not be satisfied.  In
		other words, Theorem~\ref{thm:sub-reg} gives well-posedness in
		the energy space if $\Lambda_{\infty}$ is 0 or small, while
		Theorem~\ref{thm:dgcs} provides the (DGCS)-phenomenon if
		$\Lambda_{\infty}=+\infty$.  The case where $\Lambda_{\infty}$
		is finite but large remains open.  We suspect that the
		(DGCS)-phenomenon is still possible, but our construction does
		not work.  We comment on this issue in the first part of
		section~\ref{sec:counterexamples}.
		
		Finally, Theorem~\ref{thm:dgcs} shows that assumptions
		(\ref{hp:(DGCS)-psi}) and (\ref{hp:(DGCS)-phi}) of
		Theorems~\ref{thmbibl:existence} and~\ref{thmbibl:regularity}
		are optimal also in the subcritical dissipative case with
		$\Lambda_{\infty}=+\infty$.  If initial data are in the Gevrey
		space with $\varphi(x)=x\omega(1/x)$, solutions remain in the
		same space.  If initial are in a Gevrey space corresponding to
		some $\varphi(x)\ll x\omega(1/x)$, then it may happen that for
		positive times the solution lies in the space of
		ultradistributions with $\psi(x):=x\omega(1/x)$, but not in
		the space of ultradistributions corresponding to any given
		$\psi(x)\ll x\omega(1/x)$.
	\end{em}
\end{rmk}

\setcounter{equation}{0}
\section{Heuristics}\label{sec:heuristics}

The following pictures summarize roughly the results of this paper.
In the horizontal axis we represent the time-regularity of $c(t)$.
With some abuse of notation, values $\alpha\in(0,1)$ mean that $c(t)$
is $\alpha$-H\"{o}lder continuous, $\alpha=1$ means that it is
Lipschitz continuous, $\alpha>1$ means even more regular.  In the
vertical axis we represent the space-regularity of initial data, where
the value $s$ stands for the Gevrey space of order $s$ (so that 
higher values of $s$ mean lower regularity).  The curve is
$s=(1-\alpha)^{-1}$.  

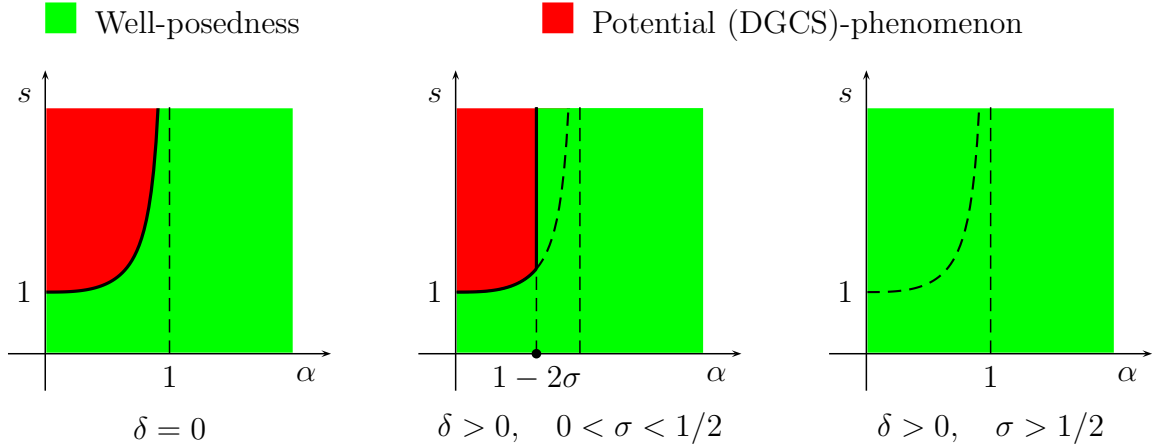
\def\grafico{\psplot{0}{0.95}{0.5 1 x 3 exp sub div}}
\psset{unit=9ex}

\noindent
\pspicture(-0.5,-1)(2.5,3)
\psclip{\psframe[linestyle=none](0,0)(2,2)}
\pscustom[fillstyle=solid,fillcolor=green,linestyle=none]{
\grafico
\psline[linestyle=none](1,2)(2,2)(2,0)(0,0)}
\pscustom[fillstyle=solid,fillcolor=red,linestyle=none]{
\grafico
\psline[linestyle=none](0,2)}
{\psset{linewidth=1.5\pslinewidth}\grafico}%
\endpsclip
\psline[linewidth=0.7\pslinewidth]{->}(-0.3,0)(2.3,0)
\psline[linewidth=0.7\pslinewidth]{->}(0,-0.3)(0,2.3)
\psline[linewidth=0.7\pslinewidth,linestyle=dashed](1,0)(1,2)
\uput[-90](2.1,0){$\alpha$}
\uput[-90](1,0){$1$}
\uput[180](0,2.1){$s$}
\uput[180](0,0.5){$1$}
\rput(1,-0.6){$\delta=0$}
\psframe*[linecolor=red](4,2.6)(4.25,2.85)
\rput[Bl](4.4,2.6){Potential (DGCS)-phenomenon}
\psframe*[linecolor=green](0,2.6)(0.25,2.85)
\rput[Bl](0.4,2.6){Well-posedness}
\endpspicture
\hfill
\pspicture(-0.5,-1)(2.5,2.5)
\psclip{\psframe[linestyle=none](0,0)(2,2)}
\pscustom[fillstyle=solid,fillcolor=green,linestyle=none]{
\grafico
\psline[linestyle=none](1,2)(2,2)(2,0)(0,0)}
\pscustom[fillstyle=solid,fillcolor=red,linestyle=none]{
\grafico
\psline[linestyle=none](0,2)}
{\psset{linewidth=1.5\pslinewidth}\grafico}%
\psframe[fillstyle=solid,fillcolor=green,linestyle=none](0.65,0)(2,2.1)
{\psset{linestyle=dashed}\grafico}%
\endpsclip
\psclip{\pscustom[linestyle=none]{\grafico\lineto(0,0.5)}}
\psline[linewidth=1.5\pslinewidth](0.65,0)(0.65,2)
\endpsclip
\psline[linewidth=0.7\pslinewidth]{->}(-0.3,0)(2.3,0)
\psline[linewidth=0.7\pslinewidth]{->}(0,-0.3)(0,2.3)
\psline[linewidth=0.7\pslinewidth,linestyle=dashed](1,0)(1,2)
\psline[linewidth=0.7\pslinewidth,linestyle=dashed](0.65,0)(0.65,2)
\uput[-90](2.1,0){$\alpha$}
\uput[-90](0.65,0){$1-2\sigma$}
\psdot(0.65,0)
\uput[180](0,2.1){$s$}
\uput[180](0,0.5){$1$}
\rput(1,-0.6){$\delta>0,\quad 0<\sigma<1/2$}
\endpspicture
\hfill
\pspicture(-0.5,-1)(2.5,2.5)
\psframe[fillstyle=solid,fillcolor=green,linestyle=none](0,0)(2,2)
\psclip{\psframe[linestyle=none](0,0)(2,2)}
{\psset{linestyle=dashed}\grafico}%
\endpsclip
\psline[linewidth=0.7\pslinewidth]{->}(-0.3,0)(2.3,0)
\psline[linewidth=0.7\pslinewidth]{->}(0,-0.3)(0,2.3)
\psline[linewidth=0.7\pslinewidth,linestyle=dashed](1,0)(1,2)
\uput[-90](2.1,0){$\alpha$}
\uput[-90](1,0){$1$}
\uput[180](0,2.1){$s$}
\uput[180](0,0.5){$1$}
\rput(1,-0.6){$\delta>0,\quad \sigma>1/2$}
\endpspicture

For $\delta=0$ we have the situation described in
Remark~\ref{rmk:dgcs-1} and Remark~\ref{rmk:dgcs-2}, namely
well-posedness provided that either $c(t)$ is Lipschitz continuous or
$c(t)$ is $\alpha$-H\"{o}lder continuous and initial data are in
Gevrey spaces of order less than or equal to $(1-\alpha)^{-1}$, and
(DGCS)-phenomenon otherwise. The same picture applies if $\delta>0$ 
and $\sigma=0$.

When $\delta>0$ and $0<\sigma< 1/2$, the full strip with
$\alpha>1-2\sigma$ falls in the well-posedness region, as stated in
Theorem~\ref{thm:sub-reg}.  The region with $\alpha<1-2\sigma$ is
divided as in the non-dissipative case.  Indeed,
Theorem~\ref{thmbibl:regularity} still provides well-posedness below
the curve and on the curve, while Theorem~\ref{thm:dgcs} provides the
(DGCS)-phenomenon above the curve.  What happens on the vertical
half-line which separates the two regions is less
clear (it is the region where $\Lambda_{\infty}$ is positive and
finite, see Remark~\ref{rmk:L-infty}).

Finally, when $\delta>0$ and $\sigma> 1/2$ well-posedness dominates
because of Theorem~\ref{thm:sup-reg}, even in the degenerate 
hyperbolic case.

Now we present a rough justification of this threshold effect.  As
already observed, existence results for problem
(\ref{pbm:eqn})--(\ref{pbm:data}) are related to estimates for
solutions to the family of ordinary differential equations
(\ref{pbm:main-ode}).

Let us consider the simplest energy function
$\mathcal{E}(t):=|\ul'(t)|^{2}+\lambda^{2}|\ul(t)|^{2}$, whose
time-derivative is
\begin{eqnarray}
	\mathcal{E}'(t) & = & -4\delta\lambda^{2\sigma}|\ul'(t)|^{2}+
	2\lambda^{2}(1-c(t))\ul(t)\ul'(t)
	\nonumber  \\
	\noalign{\vspace{0.5ex}}
	 & \leq & -4\delta\lambda^{2\sigma}|\ul'(t)|^{2}
	+\lambda(1+|c(t)|) \mathcal{E}(t).
	\label{heur:E'}
\end{eqnarray}

Since $\delta\geq 0$, a simple integration gives that
\begin{equation}
	\mathcal{E}(t)\leq \mathcal{E}(0)
	\exp\left(\lambda t+\lambda\int_{0}^{t}|c(s)|\,ds\right),
	\label{heur:E-dh}
\end{equation}
which is almost enough to establish Theorem~\ref{thmbibl:dh}.

If in addition $c(t)$ is $\omega$-continuous and satisfies the strict 
hyperbolicity condition~(\ref{hp:sh}), then (\ref{heur:E-dh}) can be 
improved to
\begin{equation}
	\mathcal{E}(t)\leq M_{1}\mathcal{E}(0)
	\exp\left(M_{2}\lambda\omega(1/\lambda)t\right)
	\label{heur:E-sh}
\end{equation}
for suitable constants $M_{1}$ and $M_{2}$.  Estimates of this kind
are the key point in the proof of both Theorem~\ref{thmbibl:existence}
and Theorem~\ref{thmbibl:regularity}.  Moreover, the (DGCS)-phenomenon
is equivalent to saying that the term $\lambda\omega(1/\lambda)$ in
(\ref{heur:E-sh}) is optimal.

Let us assume now that $\delta>0$.  If $\sigma>1/2$, or $\sigma=1/2$
and $\delta$ is large enough, then it is reasonable to expect that the
first (negative) term in the right-hand side of (\ref{heur:E'})
dominates the second one, and hence $\mathcal{E}(t)\leq
\mathcal{E}(0)$, which is enough to establish well-posedness in
Sobolev spaces.  Theorem~\ref{thm:sup-reg} confirms this intuition.

If $\sigma\leq 1/2$ and $c(t)$ is constant, then (\ref{pbm:main-ode}) 
can be explicitly integrated, obtaining that
\begin{equation}
	\mathcal{E}(t)\leq \mathcal{E}(0)\exp\left(-2\delta\lambda^{2\sigma}t\right).
	\label{heur:E-delta}
\end{equation}

If $c(t)$ is $\omega$-continuous and satisfies the strict
hyperbolicity assumption (\ref{hp:sh}), then we expect a superposition
of the effects of the coefficient, represented by (\ref{heur:E-sh}),
and the effects of the damping, represented by (\ref{heur:E-delta}).
We end up with
\begin{equation}
	\mathcal{E}(t)\leq
	M_{1}\mathcal{E}(0)\exp\left([M_{2}\lambda\omega(1/\lambda)-
	2\delta\lambda^{2\sigma}]t\right).
	\label{heur:E-conflict}
\end{equation}

Therefore, it is reasonable to expect that $\mathcal{E}(t)$ satisfies
an estimate independent of $\lambda$, which guarantees well-posedness
in Sobolev spaces, provided that
$\lambda\omega(1/\lambda)\ll\lambda^{2\sigma}$, or
$\lambda\omega(1/\lambda)\sim\lambda^{2\sigma}$ and $\delta$ is large
enough.  Theorem~\ref{thm:sub-reg} confirms this intuition.  The same
argument applies if $\sigma=0$ and $\omega(x)=Lx$, independently of
$L$ (see Remark~\ref{rmk:lip-0}).

On the contrary, in all other cases the right-hand side of
(\ref{heur:E-conflict}) diverges as $\lambda\to+\infty$, opening the
door to the (DGCS)-phenomenon.  We are able to show that it does
happen provided that $\lambda\omega(1/\lambda)\gg\lambda^{2\sigma}$.
We refer to the first part of section~\ref{sec:counterexamples} for
further comments.

\setcounter{equation}{0}
\section{Proofs of well-posedness and regularity results}\label{sec:proofs}

All proofs of our main results concerning well-posedness and
regularity rely on suitable estimates for solutions to the ordinary
differential equation (\ref{pbm:main-ode}) with initial data
\begin{equation}
	\ul(0)=u_{0},
	\hspace{3em}
	\ul'(0)=u_{1}.
	\label{eqn:ODE-data}
\end{equation}

For the sake of simplicity in the sequel we write $u(t)$ instead 
of $\ul(t)$.

\subsection{Supercritical dissipation}

Let us consider the case $\sigma\geq 1/2$. The key tool is the 
following.

\begin{lemma}\label{lemma:ODE-sup}
	Let us consider problem (\ref{pbm:main-ode})--(\ref{eqn:ODE-data})
	under the following assumptions:
	\begin{itemize}
		\item the coefficient $c:[0,+\infty)\to\re$ is measurable and
		satisfies the degenerate hyperbolicity assumption
		(\ref{hp:dh}),
	
		\item $\delta$, $\lambda$, $\sigma$ are positive real numbers
		such that
		\begin{equation}
			4\delta^{2}\lambda^{4\sigma-2}\geq\mu_{2}.
			\label{hp:super}
		\end{equation}
	\end{itemize}
	
	Then the solution $u(t)$ satisfies the following estimates.
	\begin{enumerate}
		\renewcommand{\labelenumi}{(\arabic{enumi})}
		
		\item  For every $t\geq 0$ it turns out that
		\begin{equation}
			|u(t)|^{2}\leq\frac{2}{\delta^{2}\lambda^{4\sigma}}u_{1}^{2}+
			3u_{0}^{2},
			\label{th:est-u}
		\end{equation}
		\begin{equation}
			|u'(t)|^{2}\leq
			\left(2+\frac{\mu_{2}^{2}}{\delta^{4}\lambda^{8\sigma-4}}\right)
			u_{1}^{2}+
			\frac{3\mu_{2}^{2}}{2\delta^{2}\lambda^{4\sigma-4}}u_{0}^{2}.
			\label{th:est-u'}
		\end{equation}
			
		\item Let us assume in addition that $\lambda\geq 1$ and
		$\sigma\geq 1/2$, and let $\alpha$ and $\beta$ be two real
		numbers satisfying (\ref{hp:abs}).
		
		Then for every $t\geq 0$ it turns out that
		\begin{equation}
			\lambda^{4\beta}|u'(t)|^{2}+\lambda^{4\alpha}|u(t)|^{2}\leq
			\left(2+\frac{2}{\delta^{2}}+\frac{\mu_{2}^{2}}{\delta^{4}}\right)
			\lambda^{4\beta}u_{1}^{2}+
			3\left(1+\frac{\mu_{2}^{2}}{2\delta^{2}}\right)
			\lambda^{4\alpha}u_{0}^{2}.
			\label{th:est-abs}
		\end{equation}
		
		\item In addition to the assumptions of the statement~(2), 
		let us assume also that there exists $r>0$ satisfying the 
		following three inequalities:
		\begin{equation}
			\delta\lambda^{4\sigma-2}> r\mu_{2},
			\quad\quad
			2\delta r\leq 1,
			\quad\quad
			4\delta^{2}\lambda^{4\sigma-2}\geq
			(1+2r\delta)\mu_{2}.
			\label{hp:super-r}
		\end{equation}
		
		Then for every $t\geq 0$ it turns out that
		\begin{eqnarray}
			\lambda^{4\beta}|u'(t)|^{2}+
			\lambda^{4\alpha}|u(t)|^{2} & \leq &
			\left[2\left(
			1+\frac{2\mu_{2}^{2}}{\delta^{4}}+\frac{1}{\delta^{2}}
			\right)\lambda^{4\beta}u_{1}^{2}+
			3\left(1+\frac{2\mu_{2}^{2}}{\delta^{2}}\right)
			\lambda^{4\alpha}u_{0}^{2}\right]\times
			\nonumber  \\
			 &  & \times\exp\left(-2r\lambda^{2(1-\sigma)}
			\int_{0}^{t}c(s)\,ds\right).
			\label{th:est-gevrey}
		\end{eqnarray}

	\end{enumerate}
\end{lemma}

\paragraph{\textmd{\textit{Proof}}}

Let us consider the energy $E(t)$ defined in (\ref{energy:k}). Since
$$-\frac{3}{4}|u'(t)|^{2}-\frac{4}{3}\delta^{2}\lambda^{4\sigma}|u(t)|^{2}
\leq 2\delta\lambda^{2\sigma}u(t)u'(t) \leq 
|u'(t)|^{2}+\delta^{2}\lambda^{4\sigma}|u(t)|^{2},$$
we easily deduce that
\begin{equation}
	\frac{1}{4}|u'(t)|^{2}+\frac{2}{3}\delta^{2}\lambda^{4\sigma}|u(t)|^{2}
	\leq E(t)\leq
	2|u'(t)|^{2}+3\delta^{2}\lambda^{4\sigma}|u(t)|^{2}
	\quad\quad
	\forall t\geq 0.
	\label{est:E-equiv}
\end{equation}

\subparagraph{\textmd{\textit{Statement (1)}}}

The time-derivative of $E(t)$ is
\begin{equation}
	E'(t)=-2\left(
	\delta\lambda^{2\sigma}|u'(t)|^{2}+
	\delta\lambda^{2\sigma+2}c(t)|u(t)|^{2}+
	\lambda^{2} c(t)u(t)u'(t)\right).
	\label{E'}
\end{equation}

The right-hand side is a quadratic form in $u(t)$ and $u'(t)$. The 
coefficient of $|u'(t)|^{2}$ is negative. Therefore, this quadratic 
form is less than or equal to 0 for all values of $u(t)$ and $u'(t)$ 
if and only if 
$$4\delta^{2}\lambda^{4\sigma-2}c(t)\geq c^{2}(t),$$
and this is always true because of (\ref{hp:dh}) and (\ref{hp:super}).
It follows that $E'(t)\leq 0$ for (almost) every $t\geq 0$, and hence
\begin{equation}
	\delta^{2}\lambda^{4\sigma}|u(t)|^{2}\leq E(t)\leq E(0)
	\leq 2u_{1}^{2}+3\delta^{2}\lambda^{4\sigma}u_{0}^{2},
	\label{est:E0}
\end{equation}
which is equivalent to (\ref{th:est-u}).

In order to estimate $u'(t)$, we rewrite (\ref{pbm:main-ode}) in the form
$$u''(t)+2\delta\lambda^{2\sigma}u'(t)=-\lambda^{2} c(t)u(t),$$
which we interpret as a first order linear equation with 
constant coefficients in the unknown $u'(t)$, with the 
right-hand side as a forcing term. Integrating this differential 
equation in $u'(t)$, we obtain that
\begin{equation}
	u'(t)=u_{1}\exp\left(-2\delta\lambda^{2\sigma}t\right)-
	\int_{0}^{t}\lambda^{2} c(s)u(s)
	\exp\left(-2\delta\lambda^{2\sigma}(t-s)\right)\,ds.
	\label{formula:u'}
\end{equation}

From (\ref{hp:dh}) and (\ref{th:est-u}) it follows that
\begin{eqnarray*}
	|u'(t)| & \leq & |u_{1}|+\mu_{2}\lambda^{2}\cdot
	\max_{t\in[0,T]}|u(t)|\cdot
	\int_{0}^{t}e^{-2\delta\lambda^{2\sigma}(t-s)}\,ds \\
	 & \leq &
	|u_{1}|+\frac{\mu_{2}\lambda^{2}}{2\delta\lambda^{2\sigma}}
	\left(\frac{2}{\delta^{2}\lambda^{4\sigma}}u_{1}^{2}+
	3u_{0}^{2}\right)^{1/2},
\end{eqnarray*}
and therefore
$$|u'(t)|^{2}\leq 2|u_{1}|^{2}+
\frac{\mu_{2}^{2}\lambda^{4}}{2\delta^{2}\lambda^{4\sigma}}
\left(\frac{2}{\delta^{2}\lambda^{4\sigma}}u_{1}^{2}+
3u_{0}^{2}\right),$$
which is equivalent to (\ref{th:est-u'}).

\subparagraph{\textmd{\textit{Statement (2)}}}

Exploiting (\ref{th:est-u}) and (\ref{th:est-u'}), with some simple
algebra we obtain that
\begin{eqnarray*}
	\lambda^{4\beta}|u'(t)|^{2}+\lambda^{4\alpha}|u(t)|^{2} & \leq & 
	\left(2+
	\frac{\mu_{2}^{2}}{\delta^{4}}\cdot\frac{1}{\lambda^{4(2\sigma-1)}}
	+\frac{2}{\delta^{2}}\cdot\frac{1}{\lambda^{4(\beta+\sigma-\alpha)}}
	\right)\lambda^{4\beta}u_{1}^{2}   \\
	 &  & \mbox{}+3\left(1+
	 \frac{\mu_{2}^{2}}{2\delta^{2}}\cdot
	 \frac{1}{\lambda^{4(\alpha-\beta+\sigma-1)}}
	 \right)\lambda^{4\alpha}u_{0}^{2}.
\end{eqnarray*}

All exponents of $\lambda$'s in denominators are nonnegative owing to
(\ref{hp:abs}).  Therefore, since $\lambda\geq 1$, all those fractions
can be estimated with 1.  This leads to (\ref{th:est-abs}).

\subparagraph{\textmd{\textit{Statement (3)}}}

Let us define $C(t)$ as in (\ref{defn:C(t)}).  To begin with, we prove
that in this case the function $E(t)$ satisfies the stronger
differential inequality
\begin{equation}
	E'(t)\leq -2r\lambda^{2(1-\sigma)}c(t)E(t),
	\label{est:E'}
\end{equation}
and hence
\begin{equation}
	E(t)\leq 
	E(0)\exp\left(-2r\lambda^{2(1-\sigma)}C(t)\right)
	\quad\quad
	\forall t\geq 0.
	\label{est:sup-E}
\end{equation}

Coming back to (\ref{E'}), inequality (\ref{est:E'}) is equivalent to
$$\lambda^{2\sigma}\left(\delta-r\lambda^{2-4\sigma}c(t)\right)|u'(t)|^{2}+
\delta\lambda^{2\sigma+2}(1-2r\delta)c(t)|u(t)|^{2}+
\lambda^{2}(1-2r\delta)c(t)u(t)u'(t)\geq 0.$$

As in the proof of statement~(1), we consider the whole left-hand side
as a quadratic form in $u(t)$ and $u'(t)$.  Since $c(t)\leq\mu_{2}$,
from the first inequality in~(\ref{hp:super-r}) it follows that
$$\delta\lambda^{4\sigma-2}>r\mu_{2}\geq rc(t),$$
which is equivalent to saying that the coefficient of $|u'(t)|^{2}$ is
positive. Therefore, the quadratic form is nonnegative for all values 
of $u(t)$ and $u'(t)$ if and only if
$$4\lambda^{2\sigma}\left(\delta-r\lambda^{2-4\sigma}c(t)\right)
\cdot\delta\lambda^{2\sigma+2}c(t)(1-2r\delta)\geq
\lambda^{4}c^{2}(t)(1-2r\delta)^{2},$$
hence if and only if
$$(1-2r\delta)c(t)\left[4\delta^{2}\lambda^{4\sigma-2}-
(1+2r\delta) c(t)\right]\geq 0,$$
and this follows from (\ref{hp:dh}) and from the last two inequalities
in (\ref{hp:super-r}).

Now from (\ref{est:sup-E}) it follows that
\begin{equation}
	\delta^{2}\lambda^{4\sigma}|u(t)|^{2}\leq E(t)\leq
	E(0)\exp\left(-2r\lambda^{2(1-\sigma)}C(t)\right),
	\label{est:sup-u}
\end{equation}
which provides an estimate for $|u(t)|$. In order to estimate 
$u'(t)$, we write it in the form (\ref{formula:u'}), and we estimate 
the two terms separately. The third inequality in 
(\ref{hp:super-r}) implies that $2\delta\lambda^{4\sigma-2}\geq 
r\mu_{2}$. Since $C(t)\leq\mu_{2}t$, it follows 
that
$$2\delta\lambda^{2\sigma}t\geq
r\lambda^{2-2\sigma}\mu_{2}t\geq
r\lambda^{2-2\sigma}C(t),$$
and hence
\begin{equation}
	\left|u_{1}\exp\left(-2\delta\lambda^{2\sigma}t\right)\right|\leq
	|u_{1}|\exp\left(-2\delta\lambda^{2\sigma}t\right)\leq
	|u_{1}|\exp\left(-r\lambda^{2(1-\sigma)}C(t)\right).
	\label{est:u'-1}
\end{equation}

As for the second terms in (\ref{formula:u'}), we exploit 
(\ref{est:sup-u}) and we obtain that
$$\left|\int_{0}^{t}\lambda^{2} c(s)u(s)
\exp\left(-2\delta\lambda^{2\sigma}(t-s)\right)\,ds\right| \leq
\lambda^{2}\mu_{2}\int_{0}^{t}|u(s)|
\exp\left(-2\delta\lambda^{2\sigma}(t-s)\right)\,ds$$
$$\leq \frac{\mu_{2}[E(0)]^{1/2}}{\delta\lambda^{2\sigma-2}}
\exp\left(-2\delta\lambda^{2\sigma}t\right)
\int_{0}^{t}\exp\left(-r\lambda^{2(1-\sigma)}C(s)
+2\delta\lambda^{2\sigma}s\right)ds.$$

From the first inequality in (\ref{hp:super-r}) it follows that
$$2\delta\lambda^{2\sigma}-r\lambda^{2(1-\sigma)}c(s)\geq
2\delta\lambda^{2\sigma}-r\lambda^{2(1-\sigma)}\mu_{2}\geq
\delta\lambda^{2\sigma},$$
hence
\begin{eqnarray*}
	\lefteqn{\hspace{-3em}
	\int_{0}^{t}\exp\left(-r\lambda^{2(1-\sigma)}C(s)
	+2\delta\lambda^{2\sigma}s\right)ds}  \\
	\quad\quad & \leq &
	\frac{1}{\delta\lambda^{2\sigma}}
	\int_{0}^{t}
	\left(2\delta\lambda^{2\sigma}-r\lambda^{2(1-\sigma)}c(s)\right)
	\exp\left(2\delta\lambda^{2\sigma}s-
	r\lambda^{2(1-\sigma)}C(s)\right)\,ds   \\
	 & \leq & \frac{1}{\delta\lambda^{2\sigma}}
	 \exp\left(2\delta\lambda^{2\sigma}t-r\lambda^{2(1-\sigma)}C(t)
	\right),
\end{eqnarray*}
and therefore
\begin{equation}
	\left|\int_{0}^{t}\lambda^{2} c(s)u(s)
	\exp\left(-2\delta\lambda^{2\sigma}(t-s)\right)\,ds\right|\leq
	\frac{\mu_{2}[E(0)]^{1/2}}{\delta^{2}\lambda^{4\sigma-2}}
	\exp\left(-r\lambda^{2(1-\sigma)}C(t)\right).
	\label{est:u'-2}
\end{equation}

From (\ref{formula:u'}), (\ref{est:u'-1}) and (\ref{est:u'-2}) we
deduce that
$$|u'(t)|\leq
\left(|u_{1}|+\frac{\mu_{2}[E(0)]^{1/2}}{\delta^{2}\lambda^{4\sigma-2}}\right)
\exp\left(-r\lambda^{2(1-\sigma)}C(t)\right),$$
and hence
\begin{equation}
	|u'(t)|^{2}\leq
	\left(2|u_{1}|^{2}+
	\frac{2\mu_{2}^{2}E(0)}{\delta^{4}\lambda^{8\sigma-4}}\right)
	\exp\left(-2r\lambda^{2(1-\sigma)}C(t)\right).
	\label{est:sup-u'}
\end{equation}

Finally, we estimate $E(0)$ as in (\ref{est:E0}). At this point, estimate
(\ref{th:est-gevrey}) follows from (\ref{est:sup-u'}) and
(\ref{est:sup-u}) with some simple algebra (we need to exploit that 
$\lambda\geq 1$ and assumption (\ref{hp:abs}) exactly as in the proof 
of statement~(2)).\qed

\subsubsection{Proof of Theorem~\ref{thm:sup-reg}}

Let us fix a real number $\nu\geq 1$ such that
$4\delta^{2}\nu^{4\sigma-2}\geq\mu_{2}$ (such a number
exists because of our assumptions on $\delta$ and $\sigma$).  Let us
consider the components $u_{k}(t)$ of $u(t)$ corresponding to
eigenvalues $\lk\geq\nu$.  Since $\lk\geq 1$ and
$4\delta^{2}\lk^{4\sigma-2}\geq\mu_{2}$, we can apply statement~(2) of
Lemma~\ref{lemma:ODE-sup} to these components.  If $u_{0k}$ and
$u_{1k}$ denote the corresponding components of initial data, estimate
(\ref{th:est-abs}) read as
$$\lk^{4\beta}|u_{k}'(t)|^{2}+ \lk^{4\alpha}|u_{k}(t)|^{2}\leq
\left(2+\frac{2}{\delta^{2}}+\frac{\mu_{2}^{2}}{\delta^{4}}\right)
\lk^{4\beta}|u_{1,k}|^{2}+
3\left(1+\frac{\mu_{2}^{2}}{2\delta^{2}}\right)
\lk^{4\alpha}|u_{0,k}|^{2}.$$

Summing over all $\lk\geq\nu$ we obtain exactly (\ref{th:sup-u+}).

This proves that $u_{\nu,+}(t)$ is bounded with values in
$D(A^{\alpha})$ and $u_{\nu,+}'(t)$ is bounded with values in
$D(A^{\beta})$.  The same estimate guarantees the uniform convergence
in the whole half-line $t\geq 0$ of the series defining
$A^{\alpha}u_{\nu,+}(t)$ and $A^{\beta}u_{\nu,+}'(t)$.  Since all
summands are continuous, and the convergence is uniform, the sum is
continuous as well.  Since low-frequency components
$u_{\nu,-}(t)$ and $u_{\nu,-}'(t)$ are continuous
(see Remark~\ref{rmk:low-freq}), (\ref{th:sup-reg}) is proved.\qed

\subsubsection{Proof of Theorem~\ref{thm:sup-gevrey}}

Let us fix a real number $\nu\geq 1$ such that
$4\delta^{2}\nu^{4\sigma-2}>\mu_{2}$ (such a number exists because of
our assumptions on $\delta$ and $\sigma$).  Then there exists $r>0$
such that the three inequalities in (\ref{hp:super-r}) hold true for
every $\lambda\geq\nu$.  Therefore, we can apply statement~(3) of
Lemma~\ref{lemma:ODE-sup} to all components $u_{k}(t)$ of $u(t)$
corresponding to eigenvalues $\lk\geq\nu$.  If $u_{0k}$ and $u_{1k}$
denote the corresponding components of initial data, estimate
(\ref{th:est-gevrey}) read as 
$$\left(\lk^{4\beta}|u_{k}'(t)|^{2}+
\lk^{4\alpha}|u_{k}(t)|^{2}\right)
\exp\left(2r\lk^{2(1-\sigma)} \int_{0}^{t}c(s)\,ds\right)
\leq K\left(\lk^{4\beta}|u_{1k}|^{2}+
\lk^{4\alpha}|u_{0k}|^{2}\right)
$$
for every $t\geq 0$, where $K$ is a suitable constant depending only
on $\mu_{2}$ and $\delta$.  Summing over all $\lk\geq\nu$ we obtain
exactly (\ref{th:sup-g-est}).  The continuity of $u(t)$ and $u'(t)$
with values in the suitable spaces follows from the uniform
convergence of series as in the proof of
Theorem~\ref{thm:sup-reg}.\qed

\subsection{Subcritical dissipation}

Let us consider the case $0\leq\sigma\leq 1/2$. The key tool is the 
following.

\begin{lemma}\label{lemma:ODE-sub}
	Let us consider problem (\ref{pbm:main-ode})--(\ref{eqn:ODE-data})
	under the following assumptions:
	\begin{itemize}
		\item  the coefficient $c:[0,+\infty)\to\re$ satisfies the 
		strict hyperbolicity assumption (\ref{hp:sh}) and the 
		$\omega$-continuity assumption (\ref{hp:ocont}) for some 
		continuity modulus $\omega(x)$,
	
		\item $\delta>0$, $\lambda>0$, and $\sigma\geq 0$ are real
		numbers satisfying (\ref{hp:sub-lambda}).
	\end{itemize}
	
	Then the solution $u(t)$ satisfies the following estimates.
	\begin{enumerate}
		\renewcommand{\labelenumi}{(\arabic{enumi})}
		
		\item  It turns out that
		\begin{equation}
			|u'(t)|^{2}+2\lambda^{2}\mu_{1}|u(t)|^{2}\leq
			4u_{1}^{2}+
			2\left(3\delta^{2}\lambda^{4\sigma}+
			\lambda^{2}\mu_{2}\right)u_{0}^{2}
			\quad\quad
			\forall t\geq 0.
			\label{th:sub}
		\end{equation}
	
		\item Let us assume in addition that $\lambda\geq 1$,
		$\sigma\in[0,1/2]$, and there exists a constant
		$r\in(0,\delta)$ such that
		\begin{equation}
			4(\delta-r)(\delta\mu_{1}-r\mu_{2})
			\geq
			\left[\lambda^{1-2\sigma}\omega\left(\frac{1}{\lambda}\right)\right]^{2}
			+2\delta(1+2r)
			\left[\lambda^{1-2\sigma}\omega\left(\frac{1}{\lambda}\right)\right]+
			8r\delta^{3}.
			\label{hp:sub-Lambda}
		\end{equation}
		
		Then for every $t\geq 0$ it turns out that
		\begin{equation}
			|u'(t)|^{2}+2\lambda^{2}\mu_{1}|u(t)|^{2}\leq
			\left[4u_{1}^{2}+
			2\left(3\delta^{2}\lambda^{4\sigma}+
			\lambda^{2}\mu_{2}\right)u_{0}^{2}\right]
			\exp\left(-2r\lambda^{2\sigma}t\right).
			\label{th:sub-gevrey}
		\end{equation}

	\end{enumerate}
	
\end{lemma}

\paragraph{\textmd{\textit{Proof}}}

For every $\ep>0$ we introduce the regularized coefficient
$$\cep(t):=\frac{1}{\ep}\int_{t}^{t+\ep}c(s)\,ds
\quad\quad
\forall t\geq 0.$$

It is easy to see that $\cep\in C^{1}([0,+\infty),\re)$ and satisfies 
the following estimates:
\begin{equation}
	\mu_{1}\leq\cep(t)\leq\mu_{2}
	\quad\quad
	\forall t\geq 0,
	\label{est:cep-sh}
\end{equation}
\begin{equation}
	|c(t)-\cep(t)|\leq\omega(\ep)
	\quad\quad
	\forall t\geq 0,
	\label{est:cep}
\end{equation}
\begin{equation}
	|\cep'(t)|\leq\frac{\omega(\ep)}{\ep}
	\quad\quad
	\forall t\geq 0.
	\label{est:cep'}
\end{equation}

\subparagraph{\textmd{\textit{Approximated energy}}}

For every $\ep>0$ we consider the approximated hyperbolic energy 
$\Eep(t)$ defined in (\ref{energy:Eep}).
Since
$$-\frac{1}{2}|u'(t)|^{2}-2\delta^{2}\lambda^{4\sigma}|u(t)|^{2}
\leq 2\delta\lambda^{2\sigma}u(t)u'(t) \leq 
|u'(t)|^{2}+\delta^{2}\lambda^{4\sigma}|u(t)|^{2},$$
we deduce that
\begin{equation}
	\frac{1}{2}|u'(t)|^{2}+\mu_{1}\lambda^{2}|u(t)|^{2}
	\leq\Eep(t)\leq
	2|u'(t)|^{2}+(3\delta^{2}\lambda^{4\sigma}+\lambda^{2}\mu_{2})|u(t)|^{2}
	\label{est:Eep-equiv}
\end{equation}
for every $\ep>0$ and $t\geq 0$. The time-derivative of $\Eep(t)$ is
\begin{eqnarray}
	\Eep'(t) & = & -2\delta\lambda^{2\sigma}|u'(t)|^{2}-
	2\delta\lambda^{2\sigma+2}c(t)|u(t)|^{2}
	\nonumber  \\
	\noalign{\vspace{0.5ex}}
	 &  & - 2\lambda^{2}(c(t)-\cep(t))u(t)u'(t)+
	\lambda^{2}\cep'(t)|u(t)|^{2},
	\label{Eep'}
\end{eqnarray}
hence
\begin{eqnarray}
	\Eep'(t) & \leq & -2\delta\lambda^{2\sigma}|u'(t)|^{2}-
	\left(2\delta\lambda^{2\sigma+2}c(t)-
	\lambda^{2}|\cep'(t)|\right)|u(t)|^{2}
	\nonumber  \\
	\noalign{\vspace{0.5ex}}
	 &  & +2\lambda^{2}|c(t)-\cep(t)|\cdot|u(t)|\cdot|u'(t)|.
	\label{est:Eep'}
\end{eqnarray}

\subparagraph{\textmd{\textit{Statement~(1)}}}

We claim that, for a suitable choice of $\ep$, it turns out that
\begin{equation}
	\Eep'(t)\leq 0
	\quad\quad
	\forall t\geq 0.
	\label{est:Eep'<0}
\end{equation}

If we prove this claim, then we apply (\ref{est:Eep-equiv}) with that
particular value of $\ep$ and we obtain that 
$$
\frac{1}{2}|u'(t)|^{2}+\mu_{1}\lambda^{2}|u(t)|^{2}
\leq\Eep(t)\leq\Eep(0)\leq
2u_{1}^{2}+(3\delta^{2}\lambda^{4\sigma}+\lambda^{2}\mu_{2})u_{0}^{2},$$
which is equivalent to (\ref{th:sub}).

In order to prove (\ref{est:Eep'<0}), we consider the whole right-hand
side of (\ref{est:Eep'}) as a quadratic form in $|u(t)|$ and
$|u'(t)|$.  Since the coefficient of $|u'(t)|^{2}$ is negative, this
quadratic form is less than or equal to 0 for all values of $|u(t)|$
and $|u'(t)|$ if and only if 
$$2\delta\lambda^{2\sigma}\cdot
\left(2\delta\lambda^{2\sigma+2}c(t)-\lambda^{2}|\cep'(t)|\right)
-\lambda^{4}|c(t)-\cep(t)|^{2}\geq 0,$$
hence if and only if
\begin{equation}
	4\delta^{2}\lambda^{4\sigma-2}c(t)\geq
	|c(t)-\cep(t)|^{2}+2\delta\lambda^{2\sigma-2}|\cep'(t)|.
	\label{est:quadratic}
\end{equation}

Now in the left-hand side we estimate $c(t)$ from below with
$\mu_{1}$, and we estimate from above the terms in the right-hand side
as in (\ref{est:cep}) and (\ref{est:cep'}).  We obtain that
(\ref{est:quadratic}) holds true whenever
$$4\delta^{2}\mu_{1}\geq \frac{\omega^{2}(\ep)}{\lambda^{4\sigma-2}}+
2\delta\frac{\omega(\ep)}{\lambda^{2\sigma}\ep}.$$

This condition is true when $\ep:=1/\lambda$ thanks to assumption
(\ref{hp:sub-lambda}).  This completes the proof of (\ref{th:sub}).

\subparagraph{\textmd{\textit{Statement~(2)}}}

Let us assume now that $\lambda\geq 1$ and that (\ref{hp:sub-Lambda}) 
holds true for some $r\in(0,\delta)$. In this case we claim 
that, for a suitable choice of $\ep>0$, the stronger estimate
\begin{equation}
	\Eep'(t)\leq -2r\lambda^{2\sigma}\Eep(t)
	\quad\quad
	\forall t\geq 0
	\label{est:Eep'-gevrey}
\end{equation}
holds true, hence
$$\Eep(t)\leq\Eep(0)\exp\left(-2r\lambda^{2\sigma}t\right)
\quad\quad
\forall t\geq 0.$$

Due to (\ref{est:Eep-equiv}), this is enough to deduce
(\ref{th:sub-gevrey}).  So it remains to prove
(\ref{est:Eep'-gevrey}).

Owing to (\ref{Eep'}), inequality (\ref{est:Eep'-gevrey}) is 
equivalent to
$$
2\lambda^{2\sigma}(\delta-r)|u'(t)|^{2}+
\left[2\lambda^{2\sigma+2}(\delta c(t)-r\cep(t))
-\lambda^{2}\cep'(t)-4r\delta^{2}\lambda^{6\sigma}\right]|u(t)|^{2}$$
$$
+2\left[\lambda^{2}(c(t)-\cep(t))-
2r\delta\lambda^{4\sigma}\right]u(t)u'(t)\geq 0.
$$

Keeping (\ref{hp:sh}) and (\ref{est:cep-sh}) into account, the last
inequality is proved if we show that
$$
2\lambda^{2\sigma}(\delta-r)|u'(t)|^{2}+
\left[2\lambda^{2\sigma+2}(\delta\mu_{1}-r\mu_{2})
-\lambda^{2}|\cep'(t)|-4r\delta^{2}
\lambda^{6\sigma}\right]|u(t)|^{2}$$
$$
-2\left[\lambda^{2}|c(t)-\cep(t)|+
2r\delta\lambda^{4\sigma}\right]|u(t)|\cdot|u'(t)|\geq 0.
$$

As in the proof of the first statement, we consider the whole
left-hand side as a quadratic form in $|u(t)|$ and $|u'(t)|$.  The
coefficient of $|u'(t)|$ is positive because $r<\delta$.  Therefore,
this quadratic form is nonnegative for all values of $|u(t)|$ and
$|u'(t)|$ if and only if
$$
2\lambda^{2\sigma}(\delta-r)\cdot
\left[2\lambda^{2\sigma+2}(\delta\mu_{1}-r\mu_{2})
-\lambda^{2}|\cep'(t)|-4r\delta^{2}\lambda^{6\sigma}\right]\geq
\left[\lambda^{2}|c(t)-\cep(t)|+ 2r\delta\lambda^{4\sigma}\right]^{2}.
$$

Now we rearrange the terms, and we exploit (\ref{est:cep}) and 
(\ref{est:cep'}). We obtain that the last inequality is proved if we 
show that
\begin{equation}
	4(\delta-r)(\delta\mu_{1}-r\mu_{2})\geq
	\lambda^{2-4\sigma}\omega^{2}(\ep)+
	2\delta\left(1+2r\ep\lambda^{2\sigma}\right)
	\frac{\omega(\ep)}{\ep\lambda^{2\sigma}}+
	\frac{8r\delta^{3}}{\lambda^{2-4\sigma}}.
	\label{quadr-ter}
\end{equation}

Finally, we choose $\ep:=1/\lambda$, so that (\ref{quadr-ter}) becomes
$$4(\delta-r)(\delta\mu_{1}-r\mu_{2})\geq
\left[\lambda^{1-2\sigma}\omega\left(\frac{1}{\lambda}\right)\right]^{2}+
2\delta\left(1+\frac{2r}{\lambda^{1-2\sigma}}\right)
\left[\lambda^{1-2\sigma}\omega\left(\frac{1}{\lambda}\right)\right]+
\frac{8r\delta^{3}}{\lambda^{2-4\sigma}}.$$

Since $\lambda\geq 1$ and $\sigma\leq 1/2$, this inequality follows 
from assumption (\ref{hp:sub-Lambda}).\qed

\subsubsection{Proof of Theorem~\ref{thm:sub-reg}}

Let us rewrite (\ref{defn:Linfty}) in the form
\begin{equation}
	\Lambda_{\infty}=\limsup_{\lambda\to +\infty}
	\lambda^{1-2\sigma}\omega\left(\frac{1}{\lambda}\right).
	\label{defn:Linfty-bis}
\end{equation}

Due to (\ref{hp:sub}), there exists $\nu\geq 1$ such that
(\ref{hp:sub-lambda}) holds true for every $\lambda\geq\nu$.
Therefore, we can apply statement~(1) of Lemma~\ref{lemma:ODE-sub} to
the components $u_{k}(t)$ of $u(t)$ corresponding to eigenvalues
$\lk\geq\nu$.  If $u_{0k}$ and $u_{1k}$ denote the corresponding
components of initial data, estimate (\ref{th:sub}) read as
$$|u_{k}'(t)|^{2}+2\lk^{2}\mu_{1}|u_{k}(t)|^{2}\leq 4|u_{1k}|^{2}+
2\left(3\delta^{2}\lk^{4\sigma}+ \lk^{2}\mu_{2}\right)|u_{0k}|^{2}.$$

Since $\sigma\leq 1/2$ and we chose $\nu\geq 1$, this implies that
$$|u_{k}'(t)|^{2}+2\lk^{2}\mu_{1}|u_{k}(t)|^{2}\leq 4|u_{1k}|^{2}+
2\left(3\delta^{2}+ \mu_{2}\right)\lk^{2}|u_{0k}|^{2}.$$

Summing over all $\lk\geq\nu$ we obtain exactly (\ref{th:sub-u+}).

This proves that $u_{\nu,+}(t)$ is bounded with values in $D(A^{1/2})$ and
$u_{\nu,+}'(t)$ is bounded with values in $H$.  The continuity of $u(t)$ and
$u'(t)$ with values in the same spaces follows from the uniform
convergence of series as in the proof of
Theorem~\ref{thm:sup-reg}.\qed

\subsubsection{Proof of Theorem~\ref{thm:sub-gevrey}}

Let us rewrite (\ref{defn:Linfty}) in the form
(\ref{defn:Linfty-bis}).  Due to (\ref{hp:sub}), there exists $r>0$
and $\nu\geq 1$ such that (\ref{hp:sub-Lambda}) holds true
for every $\lambda\geq\nu$.  Therefore, we can apply
statement~(2) of Lemma~\ref{lemma:ODE-sub} to the components
$u_{k}(t)$ of $u(t)$ corresponding to eigenvalues $\lk\geq\nu$.
If $u_{0k}$ and
$u_{1k}$ denote the corresponding components of initial data, estimate
(\ref{th:sub-gevrey}) reads as
$$\left(|u_{k}'(t)|^{2}+2\lk^{2}\mu_{1}|u_{k}(t)|^{2}\right)
\exp\left(2r\lk^{2\sigma}t\right)
\leq
4|u_{1k}|^{2}+ 2\left(3\delta^{2}\lk^{4\sigma}+
\lk^{2}\mu_{2}\right)|u_{0k}|^{2}.$$

Since $\sigma\leq 1/2$ and we chose $\nu\geq 1$, this implies that
$$\left(|u_{k}'(t)|^{2}+2\lk^{2}\mu_{1}|u_{k}(t)|^{2}\right)
\exp\left(2r\lk^{2\sigma}t\right)
\leq
4|u_{1k}|^{2}+ 2\left(3\delta^{2}+
\mu_{2}\right)\lk^{2}|u_{0k}|^{2}$$
for every $t\geq 0$.  Summing over all $\lk\geq\nu$ we obtain
(\ref{th:sub-g-est}) with a constant $K$ depending only on $\mu_{1}$,
$\mu_{2}$, and $\delta$.  The continuity of $u(t)$ and $u'(t)$ with
values in the suitable spaces follows from the uniform convergence of
series as in the proof of Theorem~\ref{thm:sup-reg}.\qed

\setcounter{equation}{0}
\section{The (DGCS)-phenomenon}\label{sec:counterexamples}

In this section we prove Theorem~\ref{thm:dgcs}.  Let us describe the
strategy before entering into details.  Roughly speaking, what we need
is a solution $u(t)$ whose components $u_{k}(t)$ are small at time
$t=0$ and huge at time $t>0$.  The starting point is given by the
following functions
$$b(\ep,\lambda,t):=(2\ep\lambda-\delta\lambda^{2\sigma})t-
\ep\sin(2\lambda t),$$
\begin{equation}
	w(\ep,\lambda,t):=\sin(\lambda t)\exp(b(\ep,\lambda,t)),
	\label{defn:welt}
\end{equation}
\begin{equation}
	\gamma(\ep,\lambda,t):=1+\frac{\delta^{2}}{\lambda^{2-4\sigma}}
	-16\ep^{2}\sin^{4}(\lambda t)
	-8\ep\sin(2\lambda t).
	\label{defn:gelt}
\end{equation}

With some computations it turns out that
$$w''(\ep,\lambda,t)+2\delta\lambda^{2\sigma} w'(\ep,\lambda,t)+
\lambda^{2}\gamma(\ep,\lambda,t)w(\ep,\lambda,t)=0
\quad\quad
\forall t\in\re,$$
where ``primes'' denote differentiation with respect to $t$.  As a
consequence, if we set $c(t):=\gamma(\ep,\lambda,t)$ and
$\ep:=\omega(1/\lambda)$, the function $u(t):=w(\ep,\lambda,t)$ turns
out to be a solution to (\ref{pbm:main-ode}) which grows as the
right-hand side of (\ref{heur:E-conflict}).  Unfortunately this is not
enough, because we need to realize a similar growth for countably many
components with the same coefficient $c(t)$.

To this end, we argue as in~\cite{dgcs}. We introduce a suitable 
decreasing sequence $t_{k}\to 0^{+}$, and in the interval 
$[t_{k},t_{k-1}]$ we design the coefficient $c(t)$ so that 
$u_{k}(t_{k})$ is small and $u_{k}(t_{k-1})$ is huge. Then we 
check that the piecewise defined coefficient $c(t)$ has the required 
time-regularity, and that $u_{k}(t)$ remains small for $t\in[0,t_{k}]$ 
and remains huge for $t\geq t_{k-1}$. This completes the proof.

Roughly speaking, the coefficient $c(t)$ plays on infinitely many
time-scales in order to ``activate'' countably many components, but
these countably many actions take place one by one in disjoint time
intervals.  Of course this means that the lengths $t_{k-1}-t_{k}$
of the ``activation intervals'' tend to 0 as $k\to +\infty$.  In order
to obtain enough growth, despite of the vanishing length of activation
intervals, we are forced to assume that
$\lambda\omega(1/\lambda)\gg\lambda^{2\sigma}$ as $\lambda\to
+\infty$.  In addition, components do not grow exactly as
$\exp(\lambda\omega(1/\lambda)t)$, but just more than
$\exp(\varphi(\lambda)t)$ and $\exp(\psi(\lambda)t)$.

This is the reason why this strategy does not work when
$\lambda\omega(1/\lambda)\sim\lambda^{2\sigma}$ and $\delta$ is small.
In this case one would need components growing exactly as
$\exp(\lambda\omega(1/\lambda)t)$, but this requires activation
intervals of non-vanishing length, which are thus forced to overlap.
In a certain sense, the coefficient $c(t)$ should work once again on
infinitely many time-scales, but now the countably many actions should
take place in the same time.

\subparagraph{\textmd{\textit{Definition of sequences}}}

From (\ref{hp:cex-omega}) and (\ref{hp:cex-phi-psi}) it follows that
\begin{equation}
	\lim_{x\to+\infty}x^{1-2\sigma}\omega\left(\frac{1}{x}\right)=+\infty,
	\label{hp:sigma-omega}
\end{equation}
\begin{equation}
	\lim_{x\to+\infty}
	\frac{1}{x^{1-2\sigma}\omega(1/x)}+
	\frac{\varphi(x)}{x\omega(1/x)}+
	\frac{\psi(x)}{x\omega(1/x)}=0,
	\label{hp:multifrac}
\end{equation}
and a fortiori
\begin{equation}
	\lim_{x\to+\infty}x^{1+2\sigma}\omega\left(\frac{1}{x}\right)=
	+\infty,
	\label{hp:omega-bis}
\end{equation}
\begin{equation}
	\lim_{x\to+\infty}\frac{x^{2\sigma}+\varphi(x)+\psi(x)}{x}
	\,\omega\left(\frac{1}{x}\right)=0.
	\label{hp:lim-epk}
\end{equation}

Let us consider the sequence $\{\lk\}$, which we assumed to be
unbounded.  Due to (\ref{hp:omega-bis}) and (\ref{hp:multifrac}) we
can assume, up to passing to a subsequence (not relabeled), that the 
following inequalities hold true for every $k\geq 1$:
\begin{equation}
	\lk>4\lku,
	\label{hp:l1}
\end{equation}
\begin{equation}
	\lk^{1+2\sigma}\omega\left(\frac{1}{\lk}\right)\geq
	\frac{\delta^{4}}{2^{10}\pi^{2}}\frac{1}{\lku^{2-8\sigma}}+
	\frac{4k^{2}}{\pi^{2}}\lku^{2},
	\label{hp:l22}
\end{equation}
\begin{equation}
	\lk^{1+2\sigma}\omega\left(\frac{1}{\lk}\right)\geq
	\frac{4k^{2}}{\pi^{2}}\lku^{3}
	\left(\lku^{2\sigma}+\varphi(\lku)+\psi(\lku)\right)
	\omega\left(\frac{1}{\lku}\right),
	\label{hp:l21}
\end{equation}
\begin{equation}
	\lk^{1+2\sigma}\omega\left(\frac{1}{\lk}\right)\geq
	\lku\left(\lku^{2\sigma}+\varphi(\lku)+\psi(\lku)\right)
	\omega\left(\frac{1}{\lku}\right),
	\label{hp:l23}
\end{equation}
\begin{equation}
	\frac{1}{\lk^{1-2\sigma}\omega(1/\lk)}+
	\frac{\varphi(\lk)}{\lk\omega(1/\lk)}+
	\frac{\psi(\lk)}{\lk\omega(1/\lk)}\leq
	\frac{\pi^{2}}{4k^{2}}\frac{1}{\lku^{2}}.
	\label{hp:l3}
\end{equation}

Now let us set
\begin{equation}
	t_{k}:=\frac{4\pi}{\lk},
	\hspace{3em}
	s_{k}:=\frac{\pi}{\lk}
	\left\lfloor 2\frac{\lk}{\lku}\right\rfloor,
	\label{defn:tk-sk}
\end{equation}
where $\lfloor\alpha\rfloor$ denotes the largest integer less than or 
equal to $\alpha$, and
$$\ep_{k}:=\left\{ \frac{\lk^{2\sigma}+\varphi(\lk)+\psi(\lk)}{\lk}\,
\omega\left(\frac{1}{\lk}\right)\right\}^{1/2}.$$

\subparagraph{\textmd{\textit{Properties of the sequences}}}

We collect in this section of the proof all the properties of the
sequences which are needed in the sequel.  First of all, it is clear
that $\lk\to +\infty$, hence $t_{k}\to 0$ and $\ep_{k}\to 0$ (because
of (\ref{hp:lim-epk})).  Moreover it turns out that
\begin{equation}
	\frac{t_{k-1}}{4}=\frac{\pi}{\lku}\leq s_{k}
	\leq\frac{2\pi}{\lku}=\frac{t_{k-1}}{2}.
	\label{est:sk}
\end{equation}

Keeping (\ref{hp:l1}) into account, it follows that
$$t_{k}<s_{k}<t_{k-1}
\quad\quad
\forall k\geq 1,$$
and in particular also $s_{k}\to 0$. In addition, it turns out that
\begin{equation}
	\sin(\lk t_{k})=\sin(\lk s_{k})=0
	\label{est:sin}
\end{equation}
and 
\begin{equation}
	|\cos(\lk t_{k})|=|\cos(\lk s_{k})|=1.
	\label{est:cos}
\end{equation}

Since $\sigma<1/2$, $\lk\to +\infty$, $\ep_{k}\to 0$, $t_{k}\to 0$,
keeping (\ref{hp:sigma-omega}) and (\ref{hp:multifrac}) into account,
we deduce that the following seven inequalities are satisfied provided
that $k$ is large enough:
\begin{equation}
	\frac{\delta^{2}}{\lk^{2-4\sigma}}+16\ep_{k}^{2}+8\ep_{k}
	\leq\frac{1}{2},
	\label{ineq:1}
\end{equation}
\vspace{-1\baselineskip}
\begin{equation}
	\ep_{k}\leq\frac{1}{4},
	\label{ineq:2}
\end{equation}

\begin{equation}
	16\pi\ep_{k}+\frac{16\pi\delta}{\lk^{1-2\sigma}}\leq 2\pi,
	\label{ineq:3}
\end{equation}

\begin{equation}
	\frac{1}{\lk^{1-2\sigma}\omega(1/\lk)}+
	\frac{\varphi(\lk)}{\lk\omega(1/\lk)}+
	\frac{\psi(\lk)}{\lk\omega(1/\lk)}\leq
	\frac{1}{5^{2}\cdot 2^{10}\cdot \pi^{2}},
	\label{ineq:4}
\end{equation}

\begin{equation}
	\frac{\delta^{2}}{(4\pi)^{2-4\sigma}}(2t_{k})^{1-2\sigma}
	\sup\left\{\frac{x^{1-2\sigma}}{\omega(x)}:
	x\in(0,t_{k})\right\}
	\leq\frac{1}{5},
	\label{ineq:5}
\end{equation}

\begin{equation}
	\lk^{1-2\sigma}\omega\left(\frac{1}{\lk}\right)\geq\delta^{2},
	\label{ineq:6}
\end{equation}

\begin{equation}
	\frac{2\delta^{2}}{\lku^{2-4\sigma}\omega(1/\lku)}\leq\frac{1}{5}.
	\label{ineq:7}
\end{equation}

Let $k_{0}\in\n$ be a positive integer such that (\ref{ineq:1})
through (\ref{ineq:7}) hold true for every $k\geq k_{0}$.  From
(\ref{ineq:6}) it follows that
\begin{equation}
	\ep_{k}\lk\geq\delta\lk^{2\sigma}
	\quad\quad
	\forall k\geq k_{0}.
	\label{ineq:8}
\end{equation}

From (\ref{ineq:4}) it follows that
\begin{equation}
	32\pi\frac{\ep_{k}}{\omega(1/\lk)}\leq\frac{1}{5}
	\quad\quad
	\forall k\geq k_{0}.
	\label{est:ek/omega}
\end{equation}

Since $s_{k}\geq\pi/\lku$ (see the estimate from below in 
(\ref{est:sk})), from (\ref{hp:l22}) it follows that
\begin{equation}
	\ep_{k}\lk s_{k}\geq\frac{\delta^{2}}{32}\frac{1}{\lku^{2-4\sigma}}
	\quad\quad
	\forall k\geq k_{0},
	\label{est:c'-sk-tk}
\end{equation}
\begin{equation}
	\ep_{k}\lk s_{k}\geq 2k 
	\quad\quad 
	\forall k\geq k_{0},
	\label{est:finale-3}
\end{equation}
from (\ref{hp:l21}) it follows that
\begin{equation}
	\ep_{k}\lk s_{k}\geq 2k\ep_{k-1}\lku
	\quad\quad
	 \forall k\geq k_{0},
	\label{est:finale-2}
\end{equation}
and from (\ref{hp:l3}) it follows that
\begin{equation}
	\ep_{k}\lk s_{k}\geq 
	2k\left(\lk^{2\sigma}+\varphi(\lk)+\psi(\lk)\right)
	\quad\quad
	\forall k\geq k_{0}.
	\label{est:finale-4}
\end{equation}

As a consequence of (\ref{est:finale-3}) through (\ref{est:finale-4}) it 
turns out that
\begin{equation}
	2\ep_{k}\lk s_{k}\geq k\ep_{k-1}\lku+2k\left(
	\lk^{2\sigma}+\varphi(\lk)+\psi(\lk)\right)+k
	\quad\quad
	\forall k\geq k_{0}.
	\label{est:finale-true}
\end{equation}

Finally, from (\ref{hp:l23}) it follows that
\begin{equation}
	\ep_{k}\lk\geq\ep_{k-1}\lku
	\quad\quad
	\forall k\geq k_{0}.
	\label{est:finale-1}
\end{equation}

\subparagraph{\textmd{\textit{Definition of $c(t)$ and $u(t)$}}}

For every $k\geq 1$, let $\ell_{k}:\re\to\re$ be defined by
$$\ell_{k}(t):= \frac{\delta^{2}}{t_{k-1}-s_{k}}
\left(\frac{1}{\lku^{2-4\sigma}}-\frac{1}{\lk^{2-4\sigma}}\right)
(t-s_{k})+1+\frac{\delta^{2}}{\lk^{2-4\sigma}}
\quad\quad
\forall t\in\re.$$

Thanks to (\ref{est:sin}), $\ell_{k}(t)$ is the affine function such 
that
$$\ell_{k}(s_{k})=\gamma(\ep_{k},\lk,s_{k})
\quad\quad\mbox{and}\quad\quad
\ell_{k}(t_{k-1})=\gamma(\ep_{k-1},\lku,t_{k-1}).$$

Let $k_{0}\in\n$ be such that (\ref{ineq:1}) through (\ref{ineq:7})
hold true for every $k\geq k_{0}$.  Let us set 
$$c(t):=\left\{
\begin{array}{ll}
	1 & \mbox{if $t\leq 0$},  \\
	\noalign{\vspace{0.5ex}}
	\gamma(\ep_{k},\lk,t) & \mbox{if $t\in[t_{k},s_{k}]$ for some 
	$k\geq k_{0}$},   \\
	\noalign{\vspace{0.5ex}}
	\ell_{k}(t)   & 
	\mbox{if $t\in[s_{k},t_{k-1}]$ for some $k\geq k_{0}+1$},\\
	\noalign{\vspace{0.5ex}}
	\gamma(\ep_{k_{0}},\lambda_{k_{0}},s_{k_{0}}) & 
	\mbox{if $t\geq s_{k_{0}}$.}
\end{array}
\right.$$
 
The following picture describes this definition.  The coefficient
$c(t)$ is constant for $t\leq 0$ and for $t\geq s_{k_{0}}$.  In the
intervals $[t_{k},s_{k}]$ it coincides with $\gamma(\ep_{k},\lk,t)$,
hence it oscillates, with period of order $\lambda_{k}^{-1}$ and
amplitude of order $\ep_{k}$, around a value which tends to 1.  In the
intervals $[s_{k},t_{k-1}]$ it is just the affine interpolation of the
values at the endpoints.

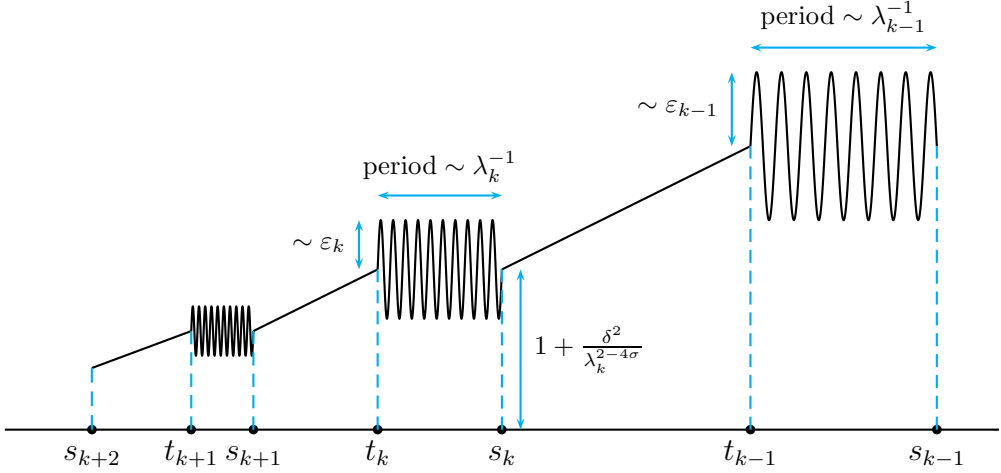
\begin{center}
	\psset{unit=9ex}
	\pspicture(-3,-2)(5,2.5)
	\psline(-2.3,-0.8)(-1.5,-0.5)
	\psplot[plotpoints=1800]{-1.5}{-1}{7200 x mul sin 0.2 mul 0.5 sub}
	\psline(-1,-0.5)(0,0)
	\psplot[plotpoints=1800]{0}{1}{3600 x mul sin 0.4 mul}
	\psline(1,0)(3,1)
	\psplot[plotpoints=1800]{3}{4.5}{1800 x mul sin 0.6 mul 1 add}
	\psline(-3,-1.3)(5,-1.3)
	\psdots(-2.3,-1.3)(-1.5,-1.3)(-1,-1.3)(0,-1.3)(1,-1.3)(3,-1.3)(4.5,-1.3)
	{\psset{linecolor=cyan,linestyle=dashed}
	\psline(-2.3,-1.3)(-2.3,-0.8)
	\psline(-1.5,-1.3)(-1.5,-0.5)
	\psline(-1,-1.3)(-1,-0.5)
	\psline(0,-1.3)(0,0)
	\psline(1,-1.3)(1,0)
	\psline(3,-1.3)(3,1)
	\psline(4.5,-1.3)(4.5,1)
	}
	\rput[B](-2.3,-1.55){$s_{k+2}$}
	\rput[B](-1.5,-1.55){$t_{k+1}$}
	\rput[B](-1,-1.55){$s_{k+1}$}
	\rput[B](0,-1.55){$t_{k}$}
	\rput[B](1,-1.55){$s_{k}$}
	\rput[B](3,-1.55){$t_{k-1}$}
	\rput[B](4.5,-1.55){$s_{k-1}$}
	\pcline[linecolor=cyan]{<->}(0,0.6)(1,0.6)
	\lput{U}{\uput[90](0,0){{\footnotesize 
	$\mbox{period}\sim\lambda_{k}^{-1}$}}}
	\pcline[linecolor=cyan]{<->}(3,1.8)(4.5,1.8)
	\lput{U}{\uput[90](0,0){{\footnotesize 
	$\mbox{period}\sim\lambda_{k-1}^{-1}$}}}
	\pcline[linecolor=cyan]{<->}(1.15,-1.3)(1.15,0)
	\lput{U}{\uput[0](0,0){{\footnotesize $1+\frac{\delta^{2}}{\lambda_{k}^{2-4\sigma}}$}}}
	\pcline[linecolor=cyan]{<->}(-0.15,0)(-0.15,0.4)
	\lput{U}{\uput[180](0,0){{\footnotesize $\sim\varepsilon_{k}$}}}
	\pcline[linecolor=cyan]{<->}(2.85,1)(2.85,1.6)
	\lput{U}{\uput[180](0,0){{\footnotesize $\sim\varepsilon_{k-1}$}}}
	\endpspicture
\end{center}

For every $k\geq k_{0}$, let us consider the solution $u_{k}(t)$ to 
the Cauchy problem
$$u_{k}''(t)+2\delta\lk^{2\sigma}u_{k}'(t)+\lk^{2}c(t)u_{k}(t)=0,$$
with ``initial'' data
\begin{equation}
	u_{k}(t_{k})=0,
	\quad\quad
	u_{k}'(t_{k})=
	\lk\exp\left((2\ep_{k}\lk-\delta\lk^{2\sigma})t_{k}\right).
	\label{uk-data}
\end{equation}

Then we set
\begin{equation}
	a_{k}:=\frac{1}{k\lk}\exp(-k\varphi(\lk)),
	\label{defn:ak}
\end{equation}
and finally
$$u(t):=\sum_{k=k_{0}}^{\infty}a_{k}u_{k}(t)e_{k}.$$

We claim that $c(t)$ satisfies (\ref{th:c-bound}) and 
(\ref{th:c-omega}), and that $u(t)$ satisfies (\ref{th:u0}) and 
(\ref{th:ut}). The rest of the proof is a verification of these claims.

\subparagraph{\textmd{\textit{Estimate and continuity of $c(t)$}}}

We prove that
\begin{equation}
	|c(t)-1|\leq\frac{1}{2}
	\quad\quad
	\forall t\geq 0,
	\label{th:c-bound-bis}
\end{equation}
which is equivalent to (\ref{th:c-bound}), and that $c(t)$ is continuous
on the whole real line.  

To this end, it is enough to check
(\ref{th:c-bound-bis}) in the intervals $[t_{k},s_{k}]$, because in 
the intervals $[s_{k},t_{k-1}]$ the function $c(t)$ is just an 
interpolation of the values at the endpoints, and it is constant for 
$t\leq 0$ and for $t\geq s_{k_{0}}$.

In the intervals $[t_{k},s_{k}]$ the function $c(t)$ coincides with
$\gamma(\ep_{k},\lk,t)$, hence from (\ref{defn:gelt}) it turns out
that
\begin{equation}
	|c(t)-1|=|\gamma(\ep_{k},\lk,t)-1|\leq
	\frac{\delta^{2}}{\lk^{2-4\sigma}}+16\ep_{k}^{2}+8\ep_{k},
	\label{th:c-bound-ter}
\end{equation}
so that (\ref{th:c-bound-bis}) follows immediately from
(\ref{ineq:1}).
	
Since the right-hand side of (\ref{th:c-bound-ter}) tends to 0 as
$k\to +\infty$, the same estimate shows also that $c(t)\to 1$ as $t\to
0^{+}$, which proves the continuity of $c(t)$ in $t=0$, the only 
point in which continuity was nontrivial.

\subparagraph{\textmd{\textit{Estimate on $c'(t)$}}}

We prove that 
\begin{equation}
	|c'(t)|\leq 32\ep_{k}\lk
	\quad\quad
	\forall t\in(t_{k},s_{k}),\ \forall k\geq k_{0},
	\label{est:c'-ts}
\end{equation}
\begin{equation}
	|c'(t)|\leq 32\ep_{k}\lk
	\quad\quad
	\forall t\in(s_{k},t_{k-1}),\ \forall k\geq k_{0}+1.
	\label{est:c'-st}
\end{equation}

Indeed in the interval $(t_{k},s_{k})$ it turns out that
$$|c'(t)|=\left|\gamma'(\ep_{k},\lk,t)\right|=
\left|-64\ep_{k}^{2}\lk\sin^{3}(\lk t)\cos(\lk t)
-16\ep_{k}\lk\cos(2\lk t)\right|$$
$$\leq 
64\ep_{k}^{2}\lk+16\ep_{k}\lk=
16\ep_{k}\lk(4\ep_{k}+1),$$
so that (\ref{est:c'-ts}) follows from (\ref{ineq:2}).

In the interval $(s_{k},t_{k-1})$ it turns out that
$$|c'(t)|=\frac{\delta^{2}}{t_{k-1}-s_{k}}
\left(\frac{1}{\lku^{2-4\sigma}}-\frac{1}{\lk^{2-4\sigma}}\right)
\leq\frac{\delta^{2}}{t_{k-1}-s_{k}}
\cdot\frac{1}{\lku^{2-4\sigma}}\leq
\frac{\delta^{2}}{s_{k}}
\cdot\frac{1}{\lku^{2-4\sigma}},$$
where the last inequality follows from the estimate from above in
(\ref{est:sk}). At this point (\ref{est:c'-st}) is equivalent  to 
(\ref{est:c'-sk-tk}).

\subparagraph{\textmd{\textit{Modulus of continuity of $c(t)$}}}

Let us prove that $c(t)$ satisfies (\ref{th:c-omega}). Since $c(t)$ 
is continuous, and constant for $t\leq 0$ and $t\geq s_{k_{0}}$, it 
is enough to verify its $\omega$-continuity in $(0,s_{k_{0}}]$. In 
turn, the $\omega$-continuity in $(0,s_{k_{0}}]$ is proved if we show 
that
\begin{equation}
	|c(t_{i})-c(t_{j})|\leq\frac{1}{5}\,\omega(|t_{i}-t_{j}|)
	\quad\quad
	\forall i\geq k_{0},\ \forall j\geq k_{0},
	\label{est:co-1}
\end{equation}
\begin{equation}
	|c(a)-c(b)|\leq\frac{1}{5}\,\omega(|a-b|)
	\quad\quad
	\forall(a,b)\in[t_{k},s_{k}]^{2},\ \forall k\geq k_{0},
	\label{est:co-2}
\end{equation}
\begin{equation}
	|c(a)-c(b)|\leq\frac{1}{5}\,\omega(|a-b|)
	\quad\quad
	\forall(a,b)\in[s_{k},t_{k-1}]^{2},\ \forall k\geq k_{0}+1.
	\label{est:co-3}
\end{equation}

Indeed, any interval $[s,t]\subseteq(0,s_{k_{0}}]$ can be decomposed
as the union of at most 5 intervals whose endpoints fit in one of the
3 possibilities above.

Let us prove (\ref{est:co-1}). From (\ref{est:sin}) it turns out that
$$|c(t_{i})-c(t_{j})|=\delta^{2}
\left|\frac{1}{\lambda_{i}^{2-4\sigma}}-\frac{1}{\lambda_{j}^{2-4\sigma}}\right|
\leq\delta^{2}
\left|\frac{1}{\lambda_{i}^{2}}-
\frac{1}{\lambda_{j}^{2}}\right|^{1-2\sigma},$$
where the inequality follows from the fact that the function $x\to 
x^{1-2\sigma}$ is $(1-2\sigma)$-H\"{o}lder continuous with constant 
equal to 1. Now from (\ref{defn:tk-sk}) it follows that
$$\delta^{2}
\left|\frac{1}{\lambda_{i}^{2}}-\frac{1}{\lambda_{j}^{2}}\right|^{1-2\sigma}
=\frac{\delta^{2}}{(4\pi)^{2-4\sigma}}|t_{i}^{2}-t_{j}^{2}|^{1-2\sigma}=
\frac{\delta^{2}}{(4\pi)^{2-4\sigma}}
|t_{i}+t_{j}|^{1-2\sigma}
\frac{|t_{i}-t_{j}|^{1-2\sigma}}{\omega(|t_{i}-t_{j}|)}
\omega(|t_{i}-t_{j}|).$$

Since $|t_{i}+t_{j}|\leq 2t_{k_{0}}$ and $|t_{i}-t_{j}|\leq 
t_{k_{0}}$, we conclude that
$$|c(t_{i})-c(t_{j})|\leq
\frac{\delta^{2}}{(4\pi)^{2-4\sigma}}(2t_{k_{0}})^{1-2\sigma}
\sup\left\{\frac{x^{1-2\sigma}}{\omega(x)}:
x\in(0,t_{k_{0}})\right\}
\omega(|t_{i}-t_{j}|),$$
so that (\ref{est:co-1}) follows from (\ref{ineq:5}).

Let us prove (\ref{est:co-2}). Since $c(t)$ is $\pi/\lk$ periodic in 
$[t_{k},s_{k}]$, for every $(a,b)\in[t_{k},s_{k}]^{2}$ there exists 
$(\oa,\ob)\in[t_{k},s_{k}]^{2}$ such that
$c(a)=c(\oa)$, $c(b)=c(\ob)$, and $|\oa-\ob|\leq\pi/\lk$. Thus from 
(\ref{est:c'-ts}) it follows that
$$|c(a)-c(b)|=|c(\oa)-c(\ob)|\leq 
32\ep_{k}\lk|\oa-\ob|=
32\ep_{k}\lk\frac{|\oa-\ob|}{\omega(|\oa-\ob|)}\omega(|\oa-\ob|),$$
so that we are left to prove that
\begin{equation}
	32\ep_{k}\lk\frac{|\oa-\ob|}{\omega(|\oa-\ob|)}\leq\frac{1}{5}.
	\label{est:co-2-bis}
\end{equation}

Due to (\ref{hp:omega-monot}), (\ref{hp:omega-monot-0}), and the fact
that $|\oa-\ob|\leq\pi/\lk$, it turns out that
$$\frac{|\oa-\ob|}{\omega(|\oa-\ob|)}\leq
\frac{\pi/\lk}{\omega(\pi/\lk)}\leq
\frac{\pi}{\lk\omega(1/\lk)},$$
so that now (\ref{est:co-2-bis}) follows from (\ref{est:ek/omega}).

Let us prove (\ref{est:co-3}). Since $c(t)$ is affine in 
$[s_{k},t_{k-1}]$, for every $a$ and $b$ in this interval it turns 
out that
$$|c(a)-c(b)|=\frac{\delta^{2}}{t_{k-1}-s_{k}}
\left(\frac{1}{\lku^{2-4\sigma}}-\frac{1}{\lk^{2-4\sigma}}\right)
|a-b|.$$

Since $s_{k}\leq t_{k-1}/2$, it follows that
$$|c(a)-c(b)|\leq
\frac{2\delta^{2}}{t_{k-1}}\frac{1}{\lku^{2-4\sigma}}\cdot|a-b|=
\frac{2\delta^{2}}{t_{k-1}}\frac{1}{\lku^{2-4\sigma}}
\cdot\frac{|a-b|}{\omega(|a-b|)}\cdot\omega(|a-b|).$$

Due to (\ref{hp:omega-monot}), (\ref{hp:omega-monot-0}), and the fact
that $|a-b|\leq t_{k-1}$, it turns out that
$$\frac{|a-b|}{\omega(|a-b|)}\leq
\frac{t_{k-1}}{\omega(t_{k-1})}\leq
\frac{t_{k-1}}{\omega(1/\lku)},$$
so that now (\ref{est:co-3}) is a simple consequence of 
(\ref{ineq:7}). 

\subparagraph{\textmd{\textit{Energy functions}}}

Let us introduce the classic energy functions
$$E_{k}(t):=|u_{k}'(t)|^{2}+\lk^{2}|u_{k}(t)|^{2},$$
$$F_{k}(t):=|u_{k}'(t)|^{2}+\lk^{2}c(t)|u_{k}(t)|^{2}.$$

Due to (\ref{th:c-bound}), they are equivalent in the sense that
$$\frac{1}{2}E_{k}(t)\leq F_{k}(t)\leq \frac{3}{2}E_{k}(t)
\quad\quad
\forall t\in\re.$$

Therefore, proving (\ref{th:u0}) is equivalent to showing that
\begin{equation}
	\sum_{k=k_{0}}^{\infty}
	a_{k}^{2}E_{k}(0)\exp(2r\varphi(\lk))<+\infty
	\quad\quad
	\forall r>0,
	\label{th:u0-equiv}
\end{equation}
while proving (\ref{th:ut}) is equivalent to showing that
\begin{equation}
	\sum_{k=k_{0}}^{\infty}
	a_{k}^{2}F_{k}(t)\exp(-2R\psi(\lk))=+\infty
	\quad\quad
	\forall R>0,\ \forall t>0.
	\label{th:ut-equiv}
\end{equation}

We are thus left to estimating $E_{k}(0)$ and $F_{k}(t)$.

\subparagraph{\textmd{\textit{Estimates in $[0,t_{k}]$}}}

We prove that 
\begin{equation}
	E_{k}(0)\leq\lk^{2}\exp(4\pi)
	\quad\quad
	\forall k\geq k_{0}.
	\label{est:ek0}
\end{equation}

To this end, we begin by estimating $E_{k}(t_{k})$.  From 
(\ref{uk-data}) we obtain that
$u_{k}(t_{k})=0$ and 
$$|u_{k}'(t_{k})|
\leq\lk\exp(2\ep_{k}\lk t_{k})=
\lk\exp(8\pi\ep_{k}),$$
so that
\begin{equation}
	E_{k}(t_{k})\leq\lk^{2}\exp(16\pi\ep_{k}).
	\label{est:ek-tk}
\end{equation}

Now the time-derivative of $E_{k}(t)$ is
$$E_{k}'(t)=-4\delta\lk^{2\sigma}|u_{k}'(t)|^{2}-
2\lk^{2}(c(t)-1)u_{k}'(t)u_{k}(t)
\quad\quad
\forall t\in\re.$$

Therefore, from (\ref{th:c-bound}) it follows that
$$E_{k}'(t)\geq -4\delta\lk^{2\sigma}E_{k}(t)
-\lk|c(t)-1|\cdot 2|u_{k}'(t)|\cdot\lk|u_{k}(t)|\geq
-\left(4\delta\lk^{2\sigma}+\frac{\lk}{2}\right)E_{k}(t)$$
for every $t\in\re$. Integrating this differential inequality in 
$[0,t_{k}]$ we deduce that
$$E_{k}(0)\leq E_{k}(t_{k})
\exp\left[\left(4\delta\lk^{2\sigma}+\frac{\lk}{2}\right)t_{k}\right].$$

Keeping (\ref{est:ek-tk}) and (\ref{defn:tk-sk}) into account, we 
conclude that
$$E_{k}(0)\leq\lk^{2}\exp\left(
16\pi\ep_{k}+\frac{16\pi\delta}{\lk^{1-2\sigma}}+2\pi\right),$$
so that (\ref{est:ek0}) follows immediately from (\ref{ineq:3}).

\subparagraph{\textmd{\textit{Estimates in $[t_{k},s_{k}]$}}}

In this interval it turns out that $u_{k}(t):=w(\ep_{k},\lk,t)$, where
$w(\ep,\lambda,t)$ is the function defined in (\ref{defn:welt}).
Keeping (\ref{est:sin}) and (\ref{est:cos}) into account, we obtain
that $u_{k}(s_{k})=0$ and
$$|u_{k}'(s_{k})|=\lk\exp(b(\ep_{k},\lk,s_{k}))=
\lk\exp\left((2\ep_{k}\lk-\delta\lk^{2\sigma})s_{k}\right).$$

Therefore, from (\ref{ineq:8}) it follows that
$$|u_{k}'(s_{k})|\geq\lk\exp(\ep_{k}\lk s_{k}),$$
and hence
\begin{equation}
	F_{k}(s_{k})=E_{k}(s_{k})\geq
	\lk^{2}\exp(2\ep_{k}\lk s_{k}).
	\label{est:fk-sk}
\end{equation}

\subparagraph{\textmd{\textit{Estimates in $[s_{k},t_{k-1}]$}}}

We prove that
\begin{equation}
	F_{k}(t_{k-1})\geq\lk^{2}
	\exp(2\ep_{k}\lk s_{k}-4\delta\lk^{2\sigma}t_{k-1}).
	\label{est:ftk-1}
\end{equation}

Indeed the time-derivative of $F_{k}(t)$ is
$$F_{k}'(t)=-4\delta\lk^{2\sigma}|u_{k}'(t)|^{2}+
\lk^{2}c'(t)|u_{k}(t)|^{2}
\quad\quad
\forall t\in(s_{k},t_{k-1}).$$

Since $c'(t)>0$ in $(s_{k},t_{k-1})$, it follows that
$$F_{k}'(t)\geq -4\delta\lk^{2\sigma}|u_{k}'(t)|^{2}\geq
-4\delta\lk^{2\sigma}F_{k}(t)
\quad\quad
\forall t\in(s_{k},t_{k-1}),$$
and hence
$$F_{k}(t_{k-1})\geq
F_{k}(s_{k})\exp\left(-4\delta\lk^{2\sigma}(t_{k-1}-s_{k})\right)\geq
F_{k}(s_{k})\exp\left(-4\delta\lk^{2\sigma}t_{k-1}\right).$$

Keeping (\ref{est:fk-sk}) into account, we have proved 
(\ref{est:ftk-1}).

\subparagraph{\textmd{\textit{Estimates in $[t_{k-1},+\infty)$}}}

We prove that
\begin{equation}
	F_{k}(t)\geq\lk^{2}\exp\left(
	2\ep_{k}\lk s_{k}-8\delta\lk^{2\sigma}t-64\ep_{k-1}\lku t\right)
	\quad\quad
	\forall t\geq t_{k-1}.
	\label{est:fk-t}
\end{equation}

To this end, let us set
$$I_{k}:=[t_{k-1},+\infty)\setminus
\bigcup_{i=k_{0}}^{k-1}\{t_{i},s_{i}\}.$$

First of all, we observe that
\begin{equation}
	|c'(t)|\leq 32\ep_{k-1}\lku
	\quad\quad
	\forall t\in I_{k}
	\label{est:c'-k-1}
\end{equation}

Indeed we know from (\ref{est:c'-ts}) and (\ref{est:c'-st}) that
$$|c'(t)|\leq 32\ep_{i}\lambda_{i}
\quad\quad
\forall t\in(t_{i},s_{i})\cup(s_{i},t_{i-1}),$$
and of course $c'(t)=0$ for every $t> s_{k_{0}}$. Now it is enough 
to observe that
$$I_{k}=(t_{k_{0}},s_{k_{0}})\cup(s_{k_{0}},+\infty)\cup
\bigcup_{i=k_{0}+1}^{k-1}[(t_{i},s_{i})\cup(s_{i},t_{i-1})],$$
and that $\ep_{i}\lambda_{i}$ is a nondecreasing sequence because of 
(\ref{est:finale-1}).

Now we observe that the function $t\to F_{k}(t)$ is continuous in 
$[t_{k-1},+\infty)$ and differentiable in $I_{k}$, with
\begin{eqnarray*}
	F_{k}'(t) & = & -4\delta\lk^{2\sigma}|u_{k}'(t)|^{2}+
	\lk^{2}c'(t)|u_{k}(t)|^{2}\\
	 & \geq & -4\delta\lk^{2\sigma}|u_{k}'(t)|^{2}-
	\frac{|c'(t)|}{c(t)}\cdot\lk^{2}c(t)|u_{k}(t)|^{2}  \\
	 & \geq & -\left(
	 4\delta\lk^{2\sigma}+\frac{|c'(t)|}{c(t)}\right)F_{k}(t).
\end{eqnarray*}

Therefore, from (\ref{est:c'-k-1}) and (\ref{th:c-bound}) it follows that
$$F_{k}'(t)\geq 
-\left(4\delta\lk^{2\sigma}+64\ep_{k-1}\lku\right)F_{k}(t)
\quad\quad
\forall t\in I_{k},$$
and hence
\begin{eqnarray*}
	F_{k}(t) & \geq & F_{k}(t_{k-1})
	\exp\left[-\left(4\delta\lk^{2\sigma}+64\ep_{k-1}\lku\right)
	(t-t_{k-1})\right]\\
	\noalign{\vspace{1ex}}
	 & \geq & F_{k}(t_{k-1})
	\exp\left[-\left(4\delta\lk^{2\sigma}+64\ep_{k-1}\lku\right)t\right]
\end{eqnarray*}
for every $t\geq t_{k-1}$.  Keeping (\ref{est:ftk-1}) into account, we
finally obtain that 
$$F_{k}(t) \geq\lk^{2} \exp\left(2\ep_{k}\lk
s_{k}-4\delta\lk^{2\sigma}t_{k-1}
-4\delta\lk^{2\sigma}t-64\ep_{k-1}\lku t\right),$$
from which (\ref{est:fk-t}) follows by simply remarking that $t\geq 
t_{k-1}$.

\subparagraph{\textmd{\textit{Conclusion}}}

We are now ready to verify (\ref{th:u0-equiv}) and
(\ref{th:ut-equiv}).  Indeed from (\ref{defn:ak}) and (\ref{est:ek0})
it turns out that
\begin{eqnarray*}
	a_{k}^{2}E_{k}(0)\exp(2r\varphi(\lk)) & \leq & 
	\frac{1}{k^{2}\lk^{2}}\exp(-2k\varphi(\lk))\cdot
	\lk^{2}\exp(4\pi)\cdot\exp(2r\varphi(\lk))    \\
	 & = & \frac{1}{k^{2}}\exp\left(4\pi+2(r-k)\varphi(\lk)\right).
\end{eqnarray*}

The argument of the exponential is less than $4\pi$ when $k$ is large
enough, and hence the series in (\ref{th:u0-equiv}) converges.

Let us consider now (\ref{th:ut-equiv}). For every $t>0$ it turns out 
that $t\geq t_{k-1}$ when $k$ is large enough. For every such $k$ we 
can apply (\ref{est:fk-t}) and obtain that
\begin{eqnarray*}
	\lefteqn{\hspace{-2em}
	a_{k}^{2}F_{k}(t)\exp(-2R\psi(\lk))} \\
	\noalign{\vspace{1ex}}
	\hspace{1em} & \geq & \frac{1}{k^{2}}\exp\left(-2k\varphi(\lk)-2R\psi(\lk)+
	2\ep_{k}\lk s_{k}-8\delta\lk^{2\sigma}t-64\ep_{k-1}\lku
	t\right).
\end{eqnarray*}

Keeping (\ref{est:finale-true}) into account, it follows that
\begin{eqnarray*}
	\lefteqn{\hspace{-2em}
	a_{k}^{2}F_{k}(t)\exp(-2R\psi(\lk))}\\
	\quad\quad & \geq &
	\frac{1}{k^{2}}\exp\left( (k-64t)\ep_{k-1}\lku+2(k-R)\psi(\lk)+
	(2k-8\delta t)\lk^{2\sigma}+k\right) \\
	 & \geq & \frac{1}{k^{2}}\exp(k)
\end{eqnarray*}
when $k$ is large enough. This proves that the series in 
(\ref{th:ut-equiv}) diverges.\qed

\label{NumeroPagine}


\begin{thebibliography}{99}

	\bibitem{CR}\textsc{G.~Chen, D.~L.~Russell}; A mathematical model
	for linear elastic systems with structural damping.  \emph{Quart.\
	Appl.\ Math.}\ \textbf{39} (1981/82), no.~4, 433--454.

	\bibitem{CT1}\textsc{S.~P.~Chen, R.~Triggiani}; Proof of
	extensions of two conjectures on structural damping for elastic
	systems.  \emph{Pacific J.\ Math.}\ \textbf{136} (1989), no.~1,
	15--55.

	\bibitem{CT2}\textsc{S.~P.~Chen, R.~Triggiani}; Characterization
	of domains of fractional powers of certain operators arising in
	elastic systems, and applications.  \emph{J.\ Differential
	Equations} \textbf{88} (1990), no.~2, 279--293.
	
	\bibitem{CT3}\textsc{S.~P.~Chen, R.~Triggiani}; Gevrey class
	semigroups arising from elastic systems with gentle dissipation:
	the case $0<\alpha<1/2$.  \emph{Proc.\ Amer.\ Math.\ Soc.}\
	\textbf{110} (1990), no.~2, 401--415.
	
	\bibitem{colombini}{\sc F.~Colombini}; Quasianalytic and
	nonquasianalytic solutions for a class of weakly hyperbolic Cauchy
	problems.  \emph{J.\ Differential Equations} \textbf{241} (2007),
	no.~2, 293--304.

	\bibitem{dgcs}{\sc F.\ Colombini, E.\ De Giorgi, S.\ Spagnolo}; Sur le
	\'{e}quations hyperboliques avec des coefficients qui ne d\'{e}pendent
	que du temp.  (French) {\em Ann.\ Scuola Norm.\ Sup.\ Pisa Cl.\ Sci.\
	(4)} \textbf{6} (1979), no.~3, 511--559.

	\bibitem{cjs}{\sc F.\ Colombini, E.\ Jannelli, S.\ Spagnolo};
	Well-posedness in the Gevrey classes of the Cauchy problem for a
	nonstrictly hyperbolic equation with coefficients depending on time.
	\emph{Ann.\ Scuola Norm.\ Sup.\ Pisa Cl.\ Sci.}\ (4) \textbf{10}
	(1983), no.~2, 291--312.

	\bibitem{gg:der-loss}{\sc M.\ Ghisi, M.\ Gobbino}; Derivative loss
	for Kirchhoff equations with non-Lipschitz nonlinear term.
	\emph{Ann.\ Scuola Norm.\ Sup.\ Pisa Cl.\ Sci.\ (5)} \textbf{8}
	(2009), no.~4, 613--646.

	\bibitem{ggh:sd}{\sc M.\ Ghisi, M.\ Gobbino, H.\ Haraux}; Local
	and global smoothing effects for some linear hyperbolic equations
	with a strong dissipation.  \emph{Trans.\ Amer.\ Math.\ Soc.}\ To
	appear.  Preprint \texttt{arXiv:1402.6595 [math.AP]}.

	\bibitem{HO}\textsc{A.~Haraux, M.~\^{O}tani}; Analyticity and
	regularity for a class of second order evolution equation.
	\emph{Evol.\ Equat.\ Contr.\ Theor.}\ \textbf{2} (2013),
	no.~1, 101--117.

	\bibitem{I1}\textsc{R.~Ikehata}; Decay estimates of solutions for
	the wave equations with strong damping terms in unbounded domains.
	\emph{Math.\ Methods Appl.\ Sci.}\ \textbf{24} (2001), no.~9,
	659--670.

	\bibitem{I2}\textsc{R.~Ikehata, M.~Natsume}; Energy decay estimates
	for wave equations with a fractional damping.  \emph{ Differential
	Integral Equations} \textbf{25} (2012), no.~9-10, 939--956.

	\bibitem{I3}\textsc{R.~Ikehata, G.~Todorova, B.~Yordanov}; Wave
	equations with strong damping in Hilbert spaces.  \emph{J.\
	Differential Equations} \textbf{254} (2013), no.~8, 3352--3368.

	\bibitem{LM}\textsc{J.-L.~Lions, E.~Magenes}, Probl\`{e}mes aux
	limites non homog\`{e}nes et applications.  Vol.  3.  (French) Travaux
	et Recherches MathŽmatiques, No.  20.  Dunod, Paris, 1970.
	
	\bibitem{MR:gevrey}\textsc{S.\ Matthes, M.\ Reissig};
	Qualitative properties of structural damped wave models.
	\emph{Eurasian Math.\ J.}\ \textbf{4}
	(2013), no.~3, 84--106.

	\bibitem{nishihara}\textsc{K.~Nishihara}; Degenerate quasilinear
	hyperbolic equation with strong damping.  \emph{Funkcial.\
	Ekvac.}\ \textbf{27} (1984), no.~1, 125--145.

	\bibitem{nishihara-decay}\textsc{K.~Nishihara}; Decay properties of
	solutions of some quasilinear hyperbolic equations with strong
	damping.  \emph{Nonlinear Anal.}\ \textbf{21} (1993), no.~1,
	17--21.

	\bibitem{reed}{\sc M.\ Reed, B.\ Simon}; \emph{Methods of Modern
	Mathematical Physics, I: Functional Analysis.  Second edition}.
	Academic Press, New York, 1980.
	
	\bibitem{shibata}\textsc{Y.~Shibata}; On the rate of decay of
	solutions to linear viscoelastic equation.  \emph{Math.\ Methods
	Appl.\ Sci.}\ \textbf{23} (2000), no.~3, 203--226.
	
\end{thebibliography}
\end{document}